\documentclass[10pt,a4paper]{amsbook}
\usepackage[latin1]{inputenc}
\usepackage{amsmath,amssymb}
\usepackage[francais]{babel}
\usepackage{color}
\usepackage{multicol}

\usepackage{latexsym,epsfig}
\usepackage{amssymb,amsmath,amsfonts}
\usepackage{exscale}
\usepackage{ifthen}
\usepackage{longtable}
\usepackage{subfigure}
\usepackage{caption}
\usepackage{cite}
\usepackage{lscape}

\newtheorem{theo}{Theorem\ }
[section]\newcommand{\fdem}{\rightline{$\Box$}}
\newtheorem{lem}[theo]{Lemma\ }

\newtheorem{prop}[theo]{Proposition\ }

\newtheorem{notation}[theo]{Notation\ }

\newtheorem{property}[theo]{Property\ }
\newtheorem{remark}[theo]{Remark.\ }

\def \H{\mathbb{H}}

\def \R{\mathbb{R}}

\def\ds{\displaystyle}

\def\fdem{$\Box$}

\begin{document}     

    \centerline{\bf \sc  \Large Convergence and Counting in  Infinite Measure}
 
\bigskip

\centerline{Fran\c{c}oise Dal'bo}
 \centerline{Marc Peign\'e \& Jean-Claude Picaud}
\centerline{Andrea Sambusetti}

\bigskip\bigskip\bigskip\bigskip

\section{Introduction}

Let $\Gamma$ be a {\em Kleinian group}. i.e. a discrete, torsionless group of isometries of a Hadamard space $X$ of negative, pinched curvature $-B^2 \leq  K_X \leq -A^2<0$, with quotient $\bar X= \Gamma \backslash X$.
%In the following $(M,g)$ is a complete connected and negatively curved manifold. Denote by $X$ it's (riemannian) universal covering with distance function $d$, by $\Gamma$ it's fundamental group acting by isometries on $X$ and fix ${\bf x},\, {\bf y}\in X$. 
This paper is concerned with two mutually related problems :  
\vspace{2mm}

 1) The description of the  distribution of the orbits of  $\Gamma$ on $X$, namely of fine asymptotic properties of the {\em orbital function} : 
$$v_{\Gamma}({\bf x},{\bf y};R):=\sharp\{\gamma\in\Gamma\slash d({\bf x},\gamma\cdot{\bf y})\leq R\}$$
for ${\bf x},{\bf y} \in X$. This has been the subject of many investigations since Margulis' \cite{Ma} (see Roblin's book \cite{R2} and   Babillot's report on \cite{BFZ}  for a clear overview). The motivations to understand the behavior of the orbital function are numerous: for instance, a simple but important invariant is its {\em exponential growth rate}
$$\delta_{\Gamma}=\limsup_{R\to\infty}{1\over R}\log(v_{\Gamma}({\bf x},{\bf y};R))$$
%(independant of ${\bf x},{\bf y}$) 
which has a major dynamical significance, since it coincides with the topological entropy of the geodesic flow when  $\bar X$ is compact, and is related to many interesting rigidity results and characterization of locally symmetric spaces, cp. \cite{ham}, \cite{bk}, \cite{bcg}.
\vspace{2mm}

   2) The pointwise behavior of the {\em Poincar\'e series} associated with $\Gamma$ :
$$P_{\Gamma}({\bf x}, {\bf y}, s):=\sum_{\gamma\in\Gamma}e^{-sd({\bf x},\gamma\cdot{\bf x})}, \qquad {\bf x}, {\bf y}\in X $$
for   and $s=\delta_{\Gamma}$, which coincides with its {\em exponent of convergence}.
%(as $P_{\Gamma}({\bf x},s)=\infty$ if $s<\delta_{\Gamma}$ and $P_{\Gamma}({\bf x},s)<\infty$ if $s>\delta_{\Gamma}$). 
The group $\Gamma$ is said to be {\em convergent} if $P_{\Gamma}({\bf x}, {\bf y}, \delta_{\Gamma})<\infty$, and {\em divergent } otherwise. Divergence can  also be understood in terms of  dynamics as, by Hopf-Tsuju-Sullivan theorem, it is equivalent to ergodicity and total conservativity of the geodesic flow with respect to the Bowen-Margulis measure   on the unit tangent bundle $U\bar X$ (see again \cite{R2} for a complete account).  
\vspace{2mm}

The regularity of the asymptotic behavior of $v_{\Gamma}$, in full generality, 
%(namely without the compact assumption) 
is well expressed in Roblin's results,
%(which involves Patterson's conformal density $(\mu_{\bf x})_{{\bf x}\in X}$ or equivalently   Bowen-Margulis  measure $m_{\Gamma}$ on the unit tangent bundle $UM$) 
which trace back to Margulis' work in the compact case:

 \begin{theo}[Margulis \cite{Ma} - Roblin \cite{R1},\cite{R2}]
 \label{MargulisRoblin}
  \label{bilan}
  ${}$

\noindent  Let $X$ be a Hadamard manifold with pinched negative curvature and $\Gamma$ a    non elementary, discrete subgroup of isometries of $X$ with {\em non-arithmetic length spectrum}\footnote{This means that the set 
$\mathcal{L} (\bar X)=\{  \ell (\gamma)\; ;\; \gamma\in\Gamma\}$  of lengths of all closed geodesics of $\bar X = \Gamma \backslash X$ is not contained in a discrete subgroup of $\R$.}: 
   
\noindent {\bf (i)} the exponential growth rate $\delta_{\Gamma}$ is a true limit;

\noindent
{\bf (ii)} if $|| m_{\Gamma}||<\infty$, then 
$v_{\Gamma}({\bf x},{\bf y};R)\sim{||\mu_{\bf x}||. ||\mu_{\bf y}||\over \delta_{\Gamma}|| m_{\Gamma}||}e^{\delta_{\Gamma}R}$;

\noindent {\bf (iii)} if $|| m_{\Gamma}||=\infty$ then $v_{\Gamma}({\bf x},{\bf y};R)=o(e^{\delta_{\Gamma}R})$,

\noindent where $(\mu_{\bf x})_{{\bf x} \in X}$ denotes the family of Patterson  conformal densities of $\Gamma$, and  $m_{\Gamma}$ the Bowen-Margulis measure on  $U\bar X$.
\end{theo}

\noindent Here,  $f\sim g$ means that $f(t)/g(t)\to 1$ when $t\to\infty$;  for $c\geq 1$, we will write $f\stackrel{c}{\asymp} g$ when   ${1\over c}\leq f(t)/g(t)\leq c$ for $t \gg 0$ (or simply $f {\asymp} g$  when the constant $c$ is not specified). The best asymptotic regularity to be expected is the existence of an equivalent, as in (ii); an explicit computation of the second term in the asymptotic development of $v_\Gamma$ is a difficult question for locally symmetric spaces (and almost 
a hopeless question in the general Riemannian setting).
\vspace{2mm}

Theorem \ref{MargulisRoblin} shows that  the key assumption for a regular behavior of $v_{\Gamma}$ is that the Bowen-Margulis measure $m_\Gamma$ is finite. This condition is clearly satisfied for {\em uniform lattices} $\Gamma$ (i.e. when $\bar X=X/\Gamma$ is compact), and  more generally for groups $\Gamma$ such that $m_{\Gamma}$ has compact support (e.g., {\em convex cocompact groups}), but it may fail for {\em nonuniform} lattices, that is  when   $\bar X=X/\Gamma$ has finite volume but is not compact.
\vspace{1mm}

 The finiteness of $m_{\Gamma}$ has a nice geometrical description in  the case of {\it geometrically finite groups}. 
 Recall that  any orbit $\Gamma\cdot{\bf x}$ accumulates on  a  closed subset  $\Lambda_{\Gamma}$ of the  geometric boundary $\partial X$ of $X$, called the {\em limit set} of $\Gamma$; the group $\Gamma $  (or the quotient manifold $\bar X$)   is said to be geometrically finite if  $\Lambda_{\Gamma}$  decomposes in the set of {\em radial limit points} (the limit points $\xi$ which are approached by orbit points in the $M$-neighborhood of any given ray   issued from $\xi$, for some  $M\!<\!\infty$) and the set of {\em bounded parabolic points} (those $\xi$  fixed  by some parabolic subgroup $P$  acting cocompactly on
  $\partial X  \setminus \{\xi\}$); for a complete  study of geometrical finiteness in variable negative curvature   see   \cite{Bow} and Proposition 1.10 in  \cite{R2}, and for a description of their topology at infinity see  \cite{DPS}. 
% As pointed out in the work of Dal'bo-Otal-Peign\'e \cite{DOP}, the geometrical finiteness  concerns actually {cuspidal ends} of the manifold $M$. 
A  finite-volume  manifold $\bar X$  is a particular case of geometrically finite manifold: it can be decomposed into a compact set and  finitely many   {\em cusps} $\bar  {\mathcal C}_i$,  
%{\color{red} et aussi d'autres bouts, qui sont comment? homeomorphes \`a des demi-espaces?} \\
i.e.  topological ends  of $\bar X$ of finite volume which are quotients  of a horoball  $\mathcal{H} _{\xi_i}$ centered at a bounded parabolic point  $\xi_i\in \partial X$ by  a maximal parabolic subgroup $P_i\subset \Gamma$ fixing $\xi_i$.  
%every cusp then lifts  on  $X$ as copies of a given fundamental domain  $C_i$ for the action of $P_i$ on   $\mathcal{H} _i$, and of all its  $\Gamma$-translates. 
\vspace{1mm}

The principle ruling the regularity of the orbital function $v_\Gamma$ of nonuniform lattices, as pointed out in \cite{DOP} and  in the following papers of the authors \cite{DPPS1}, \cite{DPPS2}, \cite{DPPS3}, is that the orbital functions 
$v_{P_i}$ (defined in a similar way as $v_{\Gamma}$)
capture the relevant information  on the wildness of the metric   inside  the cusps, which in turn may imply   $|| m_{\Gamma}||=\infty$ and the irregularity of $v_{\Gamma}$. 
%This is the way we follow along  our study. 
%One correlative principle (with regard to the dual problem) is to consider {\it convergence property} of $\Gamma$'s parabolic subgroups. Recall that to each subgroup (or even subset)   $G\subset\Gamma$ we can associate  Poincar\'e series   $P_{G}$ as well summing only over elements of $G$;   this series  converges  if $s>\delta_{G}$ and diverges if $s<\delta_{G}$ where $\displaystyle \delta_{G}=\limsup_{R\to\infty}{1\over R}\log(v_{G}({\bf x}, {\bf x}; R))$ does not depend on ${\bf x}$. Again, the group $G$ is said to be convergent if $P_{G}(\delta_{G})<\infty$ and divergent otherwise. 
In this regard, distinctive properties of  the group $\Gamma$ and of its maximal  parabolic subgroups $P_i$ are their {\em type} (convergent or divergent, as defined above) and the {\em critical gap property} (CGP), i.e. if $\delta_{P_i} < \delta_\Gamma$ for all $i$. Actually, in  \cite{DOP} it is proved that, for geometrically finite groups $\Gamma$, 
the divergence of $P_i$ implies $\delta_{P_i} < \delta_\Gamma$, and that the critical gap property implies that 
 the group $\Gamma$ is divergent  with $|| m_{\Gamma}||<\infty$. 
 On the other hand there exist geometrically finite groups $\Gamma$ which do  not satisfy the CGP: 
  we call such groups {\em exotic}, and we say that a cusp is {\em dominant} if it is associated to a parabolic subgroup $P$ with $\delta_{P} = \delta_\Gamma$.   Geometrically finite,  exotic groups may as well be convergent or divergent: in the first case, they always have $|| m_{\Gamma} || = \infty$ (by Poincar\'e recurrence and Hopf-Tsuju-Sullivan theorem, as  $|| m_{\Gamma}||<\infty$ implies total conservativity); in  the second case, in \cite{DOP} it is proved that the finiteness of  $m_{\Gamma}$ depends on the convergence of the special series 
 
$$\sum_{p\in P_i} d(x,px)e^{-\delta_\Gamma d(x,px)} < +\infty .$$

%\noindent see  again \cite{DOP}.  
%\vspace{2mm}

The main aim of this paper is to present examples of lattices $\Gamma$ for which $v_{\Gamma}$ has an irregular asymptotic behavior. According to our discussion, we will then focus  on {\em exotic, non-uniform lattices}. 
The convergence property   of exotic lattices   is an interesting question on its own: while uniform lattices (as well as convex-cocompact groups) always are divergent, the only known   examples of convergent groups, to the best of our knowledge, are given in \cite{DOP} and have infinite covolume. 
The first result of the paper is to show that both convergent and divergent  exotic lattices  do exist. 
Actually, in Section 3,  by a variation of the  construction in \cite{DOP} we obtain:

 %In the light of Theorem 1.1 part 3), this paper is devoted to provide a better understanding of that case. In %this direction we prove :
 
\begin{theo}\label{lattice convergent}
For any $N\geq 2$, there exist $N$-dimensional, finite volume manifolds of pinched negative curvature whose fundamental group $\Gamma$ is (exotic and) convergent.
\end{theo}

\noindent Constructing exotic, divergent lattices is more subtle. We prove in Section 5:

\begin{theo}\label{exotic lattice divergent}
There exist  non compact  finite area surfaces of pinched negative curvature whose group $\Gamma$ is exotic and divergent.
\end{theo}

\noindent We stress the fact that   the  examples  of Theorem \ref{lattice convergent} have infinite Bowen-Margulis measure; on the other hand,  the surfaces of 
Theorem \ref{exotic lattice divergent} can have finite or  infinite Bowen-Margulis measure, according to the chosen behaviour of the metric in the cusps.
Moreover, we believe  that the assumption on the dimension in Theorem \ref{exotic lattice divergent}  is just technical, but at present we are not able to construct similar  examples in dimension $N \geq 3$.   
\vspace{2mm}

Finally, in Section 6 we address the initial question about how irregular the orbital function can be, giving estimates for the orbital function of  a large family  of exotic lattices  with infinite Bowen-Margulis measure:

\begin{theo}\label{counting} Let $\kappa\in ]1/2, 1[$. There exist  non compact finite area surfaces with pinched negative curvature whose fundamental group $\Gamma$  satisfies the following asymptotic property: for any $ {\bf x, y} \in X$
$$  v_{\Gamma}({\bf x}, {\bf y}; R) \quad \asymp \quad  C {e^{\delta_{\Gamma}R}\over R^{1-\kappa }L(R)} \qquad {\it as} \quad R\to +\infty$$
for some slowly  varying function\footnote{A function $L(t)$ is said to be ''slowly varying'' or ''of slow growth'' if it is positive, measurable and $L(\lambda t)/L(t)\to  1$  as $t\to +\infty$ for every  $\lambda>0$.} $L: \mathbb R^+\to \mathbb R^+$ and    some constant $C=C_{\Gamma, {\bf x}, {\bf y}}>0.$    
 
\end{theo}

\noindent A  general, but more technical, result  on the orbital function of divergent, exotic lattices is given in Theorem 
\ref{countingkappa}. As far as we know, except for some sharp asymptotic formulas established by  Pollicott and Sharp  \cite{PS} for the orbital function of    normal subgroups $\Gamma$ of a cocompact Kleinian group (hence, groups which are far from being geometrically finite or with finite Bowen-Margulis measure),  these are the only  examples of such precise asymptotic behavior for the orbital function of Kleinian groups with infinite Bowen-Margulis measure.

%The consequence is  straighforward  when specifying $\Delta=1$ in the first estimation and summing $v_{\Gamma}(j+1)-v_{\Gamma}(j)$ over $j$.
%
%{\ttfamily Il parait qu'il y a un ouvrage sur les propri\'et\'es fines des fonctions\`a croissance lente, duquel on pourrait tirer l'\'equivalent dans ce qui pr\'ec\`ede.}
%

\vspace{2mm}
\noindent{\bf Remark. } This work should be considered as a companion paper to \cite{DPPS1}\& \cite{DPPS3}, where we study the asymptotic properties of the  integral version of $v_{\Gamma}$, i.e. the {\em growth function} of $X$ :
$$v_{X}({\bf x}; R):=vol_X(B(x,R)).$$
In  \cite{DPPS3}, we obtain optimal conditions on the geometry on the cusps in order that there exists a {\em Margulis function}, that is  a $\Gamma$-invariant function $c\,:\, X\to \R^+$
such that 
$$v_{X}({\bf x}; R)\sim c({\bf x})e^{\omega_X R}  \;\;\;\;\;\; \mbox{for } R \rightarrow +\infty$$
where $\omega_X$ is the exponential growth rate of the function  $v_X$ (the integral analogue of $\delta_\Gamma$).
Notice that $\omega_X$ can be different from $\delta_\Gamma$, also for lattices, as we showed in \cite{DPPS1}.

\section{Geometry of negatively curved manifolds with finite volume}

\subsection{Landscape}
\label{landscape}
 
 Additionally to those given in the introduction, we present here  notations and  familiar results about negatively curved manifolds. Amongst good references we suggest \cite{BGS},\cite{Eb},\cite{Ba} and, more specifically related to this work,   \cite{HIH} and  \cite{DPPS1}. \linebreak
In the sequel, $\bar X=X/\Gamma$ is a $N$-dimensional complete connected Riemannian manifold  with metric $g$ whose sectional curvatures satisfy : $-B^2\leq K_X\leq -A^2<0$ for fixed constants $A$ and $B$;    we will assume 
$0<A\leq 1\leq B$ since in most examples $g$ will be obtained by perturbation of a  hyperbolic one and the curvature will  equal  $-1$ on large subsets of $\bar X$. 
\vspace{2mm}

The family of normalized  distance functions : 
  $$d({\bf x}_0,{\bf x})-d({\bf x}, \cdot )$$   converges uniformly on compacts to the  {\em Busemann function}  $ \mathcal B_{\xi}( {\bf x}_0, \cdot ) $  
for  ${\bf x}\to\xi\in \partial X$.  \linebreak
The {\em horoballs} $\mathcal{H}_{\xi}$       (resp. the  {\em horospheres} $\partial \mathcal{H}_{\xi}$)  centered at $\xi$  are   the sup-level sets (resp. the  level sets) of the function $ \mathcal B_{\xi}( {\bf x}_0, \cdot ) $.
 Given a  horosphere $\partial \mathcal{H}_{\xi}$ passing through a point ${\bf x}$,  we also set,  for all $t\in\R$,  
  $$\mathcal{H}_{\xi}(t):=\{{\bf y}\slash  \mathcal B_{\xi}( {\bf x}_0, {\bf y} )   \geq \mathcal B_{\xi}({\bf x}_0, {\bf x})+t \}$$
   resp. $\partial\mathcal{H}_{\xi}(t):=\{{\bf y}\slash \mathcal B_{\xi}( {\bf x}_0, {\bf y} ) = \mathcal B_{\xi}( {\bf x}_0, {\bf x} )  + t\}$.  
%and simply set $\mathcal{H}=  \mathcal{H}(0)$, $\partial\mathcal{H}= \partial\mathcal{H}(0)$; 
We will   refer to $t = \mathcal B_{\xi}( {\bf x}_0, {\bf y}) -  \mathcal B_{\xi}( {\bf x}_0,{\bf x})$ as to the {\em height} of ${\bf y}$  (or of the horosphere  $\partial\mathcal{H}_{\xi} (t)$) in $\mathcal{H}_{\xi}$. Also, when no confusion is possible, we will drop   the index $\xi \in \partial X$ denoting the center of the horoball. 
Also, recall that the Busemann function satisfies the fundamental cocycle relation  
$$\mathcal B_x(  {\bf x}, {\bf z}) = B_x(  {\bf x}, {\bf y}) + B_x(  {\bf y}, {\bf z})$$
which will be crucial in the following.

\vspace{1mm}
 
An origin ${\bf x}_0 \in X$ being fixed,  the Gromov product between    $x, y \in \partial X \cong \mathbb S^{n-1}$, \linebreak $x \neq y$, is  defined as 
$$(x\vert y)_{{\bf x}_0} = {\mathcal B_x({\bf x}_0, {\bf z})+\mathcal B_y({\bf x}_0, {\bf z})\over 2}$$
 where ${\bf z}$ is any point on the geodesic $(x, y)$ joining $x$ to $y$ ; 
 then, for any    $ 0< \kappa^2\leq A^2$  the expression
$$ D(x, y)= e^{ -\kappa(x\vert y)_{{\bf x}_0} }$$
defines  a distance on $\partial X $, cp  \cite{bour}. Recall that for any $\gamma \in \Gamma$ one gets 
\begin{equation} \label{TAF1}
D(\gamma  \cdot x), \gamma \cdot y)= e^{-{\kappa\over 2}  \mathcal B_x(\gamma^{-1}\cdot {\bf x}_0, {\bf x}_0)} e^{-{\kappa\over 2}  \mathcal B_y(\gamma^{-1}\cdot {\bf x}_0, {\bf x}_0)} D(x, y).
\end{equation}
In other words, the isometry $\gamma$ acts on $(\partial X, D)$ as a conformal transformation with coefficient of conformality 
$ 
\vert \gamma'(x)\vert = e^{-\kappa  \mathcal B_x(\gamma^{-1}\cdot {\bf x}_0, {\bf x}_0)}
$ 
at  $x$, since equality (\ref{TAF1}) may be rewritten
\begin{equation} \label{TAF2}
D(\gamma  \cdot x), \gamma \cdot y)= \sqrt{\vert \gamma'(x)\vert \vert \gamma'(y)\vert } \, D(x, y).
\end{equation}
 %  ${\mathfrak b} (\gamma, x):= \mathcal B_x(\gamma^{-1}\cdot {\bf x}_0, {\bf x}_0)$ satisfies the cocycle relation ${\mathfrak b} (\gamma_1\gamma_2, x)= {\mathfrak b} (\gamma_1, \gamma_2 \cdot x))+{\mathfrak b} (\gamma_2, x)$. \\
%\vspace{1mm}

 Recall that $\Gamma$ is a torsion free nonuniform lattice  acting on $X$ by hyperbolic or parabolic isometries.
For any $\xi \in \partial X$, denote by 
$(\psi_{\xi, t})_{t\geq 0} $ the radial semi-flow  defined as follows : for any ${\bf x} \in X$,
the point  $\psi_{\xi, t}({\bf x})$  lies on the geodesic ray $[{\bf x}, \xi)$  at distance $t$ from ${\bf x}$. \linebreak
%%For any horosphere $partial \mathcal{H}$  centered at $\xi$, we set $\partial \mathcal{H}(t)= \psi_{\xi,  t}({\partial \mathcal H})$ ;  furthermore,  if $A$ is a subset of $\partial \mathcal{H} $ and  $I$ and interval  in $\R$, we denote by $A\times I$ the set 
% $\displaystyle{ \bigcup_{s\in I}\ \psi_{\xi, s}(A)}$ (with the obvious convention $A\times I=\emptyset$ when  either $F =\emptyset$ or $  I=\emptyset$).
% 
By  classical comparison theorems on Jacobi fields (see for instance \cite{HIH}),   the differential of $\psi_{\xi, t} : \partial {\mathcal H}_{\xi}  \to \partial {\mathcal H}_{\xi} (t)$  satisfies
$
e^{-Bt}||v||\leq ||d\psi_{\xi, t}(v)||\leq e^{-At}||v||
$
 for any   $t\geq 0$  and any vector $v$ in the tangent space $ T(\partial \mathcal{H}_{\xi})$;  
consequently,  if $\mu_t$ is the Riemannian measure induced on $\partial \mathcal{H}_{\xi}(t)$ by the metric of $X$, we have, for any Borel set  $F\subset \partial \mathcal{H}_{\xi}$  
$$
e^{-B(N-1)t}\mu_0(F)\leq \mu_t(\psi_{\xi,  t}(F))= \int _{F} |Jac(\psi_{\xi,  t})(x)|d\mu_0(x) \leq e^{-A(N-1)t}\mu_0(F).
$$

As  $\bar X= \Gamma \backslash X$ is non compact and  $vol(\bar X)<\infty$,  the manifold $\bar X$ can be decomposed into a disjoint union of a relatively compact subset  $\bar {\mathcal K} $ and finitely many cusps $\bar {\mathcal C}_1, ..., \bar {\mathcal C}_l$,  each of which is a quotient   of   a horoball   $\mathcal{H}_{\xi_i}$, centered at some boundary point $\xi_i$,   by a  maximal parabolic subgroup $P_i$. 
 As a consequence of Margulis' lemma, we can choose  the family $(\mathcal{H}_{\xi_i})_{1\leq i\leq l}$ so that any two $\Gamma$-translates of the $\mathcal{H}_{\xi_i}$ are either disjoint or coincide (cp. \cite{R2}, Proposition 1.10); we call these $\mathcal{H}_{\xi_i}$ a {\em fundamental system of horospheres} for $X$. 
 Accordingly, the Dirichlet  domain ${\mathcal D}$ of $\Gamma$   centered at the base point ${\bf x}_0$  can be decomposed into a disjoint union    ${\mathcal D} =  {\mathcal K}  \cup {\mathcal C}_1 \cup \cdots \cup {\mathcal C}_l$, where 
$ {\mathcal K}$ is a convex, relatively compact set containing ${\bf x}_0$ in its interior and projecting to $\bar {\mathcal K} $, and  the ${\mathcal C}_i$   are connected fundamental domains for the action of $P_i$ on $\mathcal{H}_{\xi_i}$, projecting to    $\bar {\mathcal C}_i$. We let  ${\mathcal S}_i = {\mathcal D} \cap \partial \mathcal{H}_{\xi_i}$ be the corresponding, relatively compact fundamental domain for the action of $P_i$ on $\partial \mathcal{H}_{\xi_i}$, so that 
 ${\mathcal C}_i = {\mathcal D} \cap  \mathcal{H}_{\xi_i} \simeq   {\mathcal S}_i \times \mathbb{R}_+$.
\vspace{1mm}

Fixing  an end $\bar{\mathcal{C}}$, and omitting in what follows the index $i$, let  $\mu_t$ be  the Riemannian measure induced by the Riemannian metric on    the horosphere $\partial\mathcal{H}_{\xi}(t)$ corresponding to  $\bar{ \mathcal C}$. In \cite{DPPS1} we defined the {\it horospherical area function} associated with the cusp $\bar{\mathcal{C}}$ as :
$$\mathcal{A}(t)=\mu_t(P \backslash\partial\mathcal{H}_{\xi}(t)  )=\mu_t(\psi_{\xi, t}(\mathcal{S})).$$
This function depends on the choice  of the initial horosphere $\partial\mathcal{H}_{\xi}$ for the end $\bar{\mathcal{C}}$, and the following result shows that this dependance is unessential for our counting problem :
\begin{prop} \label{aireetcroissancedeP}{\bf \cite{DPPS1}} There exists a constant $c=c(A,B,diam({\mathcal S}))$ such that
$$v_P(R)\stackrel{c}{\asymp}{1\over {\mathcal A}({R\over 2})}.$$
\end{prop}
\noindent This weak equivalence is the key to relate the  irregularity of the metric in the cusp to the irregular asymptotic  behaviour of the orbital function of $P$.  The second crucial step of  our work  will then be to describe precisely  the contribution of $v_P$ in the  asymptotic behavior of $v_{\Gamma}$ assuming $\delta_P=\delta_{\Gamma}$.

\subsection{ Cuspidal geometry}
\label{cuspidalgeometry}

The strategy to construct  examples with  irregular orbital  functions  as in Theorems \ref{lattice convergent}, \ref{exotic lattice divergent} and $\ref{counting}$   is to perturb in a suitable manner the metric  of a finite volume hyperbolic manifold   $  \Gamma  \backslash \mathbb{H}^N $  in a cuspidal end $  P \backslash {\mathcal H} $. 
If  ${\mathcal H}   = \{{\bf y}\slash \mathcal B_{\xi}( {\bf x}_0, {\bf y} ) \geq t_0 \}$,   
the hyperbolic metric writes    on 
${\mathcal H}    \simeq  \partial {\mathcal H}    \times \mathbb{R}_+ \equiv \R^N  \times \mathbb{R}_+$  as 
 $g=e^{- (t-t_0)}d{\bf x}^2+dt^2$   in the horospherical coordinates
  ${\bf y}= ({\bf x},t)$,
   where $d{\bf x}^2$ denotes the induced flat Riemannian metric of $\partial {\mathcal H}  $  
%for ${\bf x} \in \partial {\mathcal H}_\xi$ 
and $t = \mathcal B_{\xi}( {\bf x}_0, {\bf y} ) -t_0$.
We will consider a new metric $g$ in $P \backslash {\mathcal H}_\xi$  whose lift to ${\mathcal H} $ writes,  
in the same coordinates,  as
$$g=\tau^2(t)d{\bf x}^2+dt^2.$$

\noindent  We extend this metric by $\Gamma$-invariance to $\Gamma  {\mathcal H} $ and produce  a new Hadamard  space $(X,g) $ with quotient $\bar X = \Gamma  \backslash X$.
The new manifold $\bar X$ has again finite volume, provided that  
$$\int_0^{+\infty} \tau^{n-1}(t) dt < \infty$$

\noindent and the end  $\bar {\mathcal C} = P \backslash {\mathcal H} $ is a new cusp; we call   the function $\tau$ the {\em analytic profile} of the cusp  $\bar {\mathcal C} $.   \\

\vspace{-3mm}
\noindent The horospherical area function ${\mathcal A}$ associated with the profile $\tau$ satisfies ${\mathcal A}\asymp \tau^{N-1}$;    by Proposition \ref{aireetcroissancedeP}, it implies that :
\vspace{1mm}

\noindent (a)  the parabolic group  $P$  has critical exponent $\delta_P={(N-1)\omega_\tau\over 2} $, for 
$\displaystyle  \omega_\tau:=\limsup_{R\to +\infty}{1\over R}|\ln(\tau(R))|,$

\noindent (b)    $P$    is convergent if and only if
$\displaystyle \int_{0}^{+\infty}{e^{-\omega_\tau(N-1)t}\over \tau^{N-1}(t)}{\rm d}t<\infty.$

\noindent Also, notice that the sectional curvatures at   $({\bf x},t) $ are given by $K_{({\bf x},t)}({\partial\over \partial x_i},{\partial\over \partial x_j})=-\left({\tau'\over \tau}\right)^2$ and $K_{({\bf x},t)}({\partial\over \partial x_i},{\partial\over \partial t})=-{\tau''\over \tau}$
(see \cite{BC}).  
 \vspace{2mm}

In sections 4, 5  and 6 we will apply this strategy to specific analytic profiles $\tau$ associated with a function $\tau$, and depending on additional real parameters $a,b,\eta$,  defined as follows. 

\noindent For any convex function  $\tau\,:\,\R\to\R^+$ with $\int_0^{+\infty} \tau^{n-1}(t) dt < \infty$ and satisfying the conditions
 \vspace{-3mm}
 
\begin{equation}
\label{conditionsur1}
  \forall\; t\leq t_0 \quad \tau(t)=e^{- (t-t_0)} 
\end{equation}  

\vspace{-5mm}

\begin{equation}
\label{conditionsur2}
 A<\tau''/\tau<B 
\end{equation}  

\vspace{-4mm}

\begin{equation}
\label{conditionsur3}
 \omega_{\tau} = \limsup_{t\to +\infty} \frac{|\ln(\tau(t))|}{t}<B 
 \end{equation}

 \noindent we will set, $a\geq t_0$
     $$\tau_{a}(t)=e^{-a}\tau(t-a) \;\;\;\; \mbox{ for } t\in\R.$$ This profile  defines a manifold $\bar X$ with a cusp $\bar {\mathcal C}$ which is hyperbolic at height  less than $a$ and then   has (renormalized)  profile equal to $\tau$.\\
  Moreover,  given  parameters $a \geq t_0$, $b\geq 0$ and   $\eta\in ]0,A[$, a straightforward calculus proves the existence of  a profile $\tau_{a,b,\eta}$ such that 
$$\left\{\begin{array}{ll}A-\eta<  \tau_{ a,b,\eta}''(t) / \tau_{a,b,\eta} (t) <B+\eta&  \mbox{for all } t,  \vspace{1mm} \\
                                    \tau_{ a,b,\eta}(t)=e^{-t}\quad & {\rm for} \quad t\leq a,  \vspace{1mm}  \\
                                   \tau_{ a,b,\eta}(t)=e^{-\omega_{\tau}t}\quad & {\rm for}\quad t\in [\Delta+a,\Delta+a+b],  \vspace{1mm} \\
                                   \tau_{ a,b,\eta}(t)=e^{-(2\Delta+a+b)}\tau(t-(2\Delta+a+b))\quad & {\rm for}\quad t\geq 2\Delta+a+b.\end{array}\right.$$
for some constant $\Delta=\Delta(A,B,\eta) \gg0$.  The function $\tau_{a,b, \eta}$ defines  a  cusp $\bar {\mathcal C}$  which is  hyperbolic till height   $a$, with   constant curvature $-\omega_{\tau}^2$ in a band of width $b$ at height $\Delta+a$, and which then has asymptotic  (suitable renormalized) profile $\tau$  after height  $2\Delta+a+b$.

\section{Proof of Theorem \ref{lattice convergent} :   construction of convergent lattices}

An example of  manifold of negative curvature with infinite volume and whose fundamental group is convergent is due to Dal'Bo-Otal-Peign\'e \cite{DOP}. We propose in this section a variation of their  argument to produce a convergent  nonuniform  {\em lattice}. 
We will consider a  finite volume  hyperbolic manifold $\Gamma \backslash \mathbb H^N$ with   {\em one} cuspidal end $ P \backslash  \mathcal{H}$ and  deform   the metric in this end as explained before  to obtain a new Hadamard space   $(X,g)$ such that the quotient $\bar X=\Gamma\backslash X$ has finite volume  and a dominant  cusp $\bar{\mathcal{C}} = P \backslash  \mathcal{H}$ with a {\em convergent} parabolic group $P$,   whose exponent $\delta_P$ :

\noindent - is greater than the Poincar\'e exponent of $\Gamma$ acting on  $\mathbb H^N$, that is $N-1$; 

\noindent - equals the Poincar\'e exponent  $\delta_\Gamma$ of $\Gamma$ corresponding to the new metric.

\noindent For this, we choose $\tau$ satisfying the conditions
(\ref{conditionsur1}), (\ref{conditionsur2}), (\ref{conditionsur3}) 
with  $ \delta_P={(N-1)\omega_\tau \over 2}>(N-1)$ and 
$\int_0^{+\infty}{e^{-\omega_\tau (N-1)t}\over  \tau (t) ^{N-1}}{\rm d}t<\infty$, 
and we consider
the profile $\tau_{ a}$ for some $a>0$ to be precised.
\linebreak
%With the same notations as in \ref{landscape} for $\mathcal{K}$ and the fundamental horosphere $ \mathcal{H}_0$ of $\bar X$,  let $a>0$ and set $c=2 \mbox{diam} (  \mathcal{K})$; then, we replace the hyperbolic metric on $\bar{\mathcal{C}}_0(a+c) := P \backslash  \mathcal{H}_{0} (a+c) $ by a  metric $g$ whose  profile $T$ satisfies $ \delta_P={(N-1)\omega_T\over 2}>(N-1)$ and $\int_0^{+\infty}{e^{-(N-1)\omega_T t}\over  T(t) ^{N-1}}{\rm d}t<\infty$. 
Remark  that  the first condition can be satisfied only if $\omega_\tau>2$ which requires $B^2/A^2>4$.  
We will denote by $d$ the distance on $(X,g)$  corresponding to this new metric, and by $d_0$ the hyperbolic distance.
We emphasize that  the perturbation of the metric will not change neither the  algebraic structure of the  groups $\Gamma$ and $P$, nor the horospheres $\mathcal{H}(t)$ (which are only modified in size and not   as subsets of $\H^n$)  and  their radial flow; however,    the orbital function has different behavior before and after perturbation. 
% we set $P_0$ for the parabolic subgroup associated with the cusp in the initial hyperbolic framework and $P$ after modifying the profile of the cusp. 
 \vspace{2mm}

Now, we need to introduce a natural decomposition of geodesic segments according to their excursions in the cusp, which will enable us to  encode elements of $\Gamma$ by sequences of parabolic elements  travelling far in the cusp and elements   staying in a fixed, compact subset.  We use the  same notations as in \ref{landscape} for the compact subset $\mathcal{K}$, the fundamental horosphere $ \mathcal{H} =\mathcal{H}_{\xi_1}  $ of $X$, and  the Borel fundamental domain  $\mathcal{S}$ be for the action of $P$ on $\partial \mathcal{H}$. \linebreak
% and let  $c=2 \mbox{diam} (  \mathcal{K})$.
%For $a>0$ let $\mathcal{C}_a\subset \mathcal{C}$ the set of points   of $\mathcal C$  at distance $a$ from $\partial\mathcal{C}$ for the hyperbolic metric. In the sequel, we encode elements of $\Gamma $ (and those of $\Gamma$ as well)  using translates of $\mathcal{H}(a)$. 
%For, we recall that ${\bf x}_0\in \partial\mathcal{C}$ is fixed. F
Let $h>0$: for every $\gamma\in\Gamma$, the geodesic segment $[{\bf x}_0,\gamma\cdot{\bf x}_0]$  intersects $r=r(\gamma)$ disjoint translates $g.\mathcal{H}(h)$ (with the convention $r=0$ if the intersection with $\bigcup_{g\in\Gamma}g.\mathcal{H}(h)$ is empty). In case $r\geq 1$, denote by ${\bf z}_1^-,\cdots,{\bf z}_r^-$ (resp. ${\bf z}_1^+,\cdots,{\bf z}_r^+$) the hitting (resp. exit) points of the oriented geodesic segment $[{\bf x}_0,\gamma\cdot{\bf x}_0]$ with translates of $\mathcal{H}(h)$ in this order.  Hence we get
\small
\vspace{-2mm}

$$
[{\bf x}_0, \gamma\cdot{\bf x}_0] \cap \left(\bigcup_{g \in \Gamma} g\cdot \mathcal H (h)\right)= [{\bf z}_1^-, {\bf z_1}^+] \cup \cdots \cup [{\bf z}_{r}^-, {\bf z}_{r}^+].
$$

\normalsize
%Those points are defined as there is no orbit points of ${\bf x}_0$ in $\bigcup_{g\in\Gamma_0}g.\mathcal{H}(a)$. 
%When $r\geq 1$,   we can define in a corresponding manner points 
% ${\bf y}_1^-, \cdots,{\bf y}_r^-$ (resp. ${\bf y}_1^+, \cdots,{\bf y}_r^+$) hitting (resp. exit) points of iterates of $\mathcal{H}$. 
\noindent Accordingly, when $r\geq 1$ we can define the points   ${\bf y}_1^-, {\bf y}_1^+ \cdots {\bf y}_{r}^-, {\bf y}_{r}^+$  on  $[{\bf x_0}, g\cdot {\bf x}_0]$  such that, for any $1\leq i\leq r$, the geodesic segment $[{\bf y}_{i}^-, {\bf y}_{i}^+] $ is the connected component of $[{\bf x_0}, \gamma \cdot {\bf x}_0]\cap \left(\bigcup_{g \in \Gamma} g\cdot \mathcal H\right) $  containing $[{\bf z}_{i}^-, {\bf z}_{i}^+]$. We also set ${\bf y}_0^+={\bf x}_0$ and ${\bf y}_{r+1}^-=\gamma\cdot{\bf x}_0$.
 
\noindent  With these notations, there exist uniquely determined isometries $g_1,\cdots,g_r\in\Gamma$ and \linebreak $p_1,\cdots, p_r\in P$ such that  
$ {\bf y}_1^-\in g_1 \cdot \mathcal{S} $, $ {\bf y}_1^+ \in g_1p_1 \cdot \mathcal{S} , \cdots, 
 {\bf y}_r^+ \in g_1p_1\cdots g_rp_r \cdot \mathcal{S} $.  Finally, we define   $g_{r+1}$ by the relation  
$$\gamma=g_1p_1\cdots g_rp_rg_{r+1}$$

\noindent which we call    the {\em horospherical decomposition  of $\gamma$ at height $h$}. Notice that this decomposition  depends only on the initial hyperbolic metric.  
We also set ${\bf x}_0^+={\bf x}_0, {\bf x}_{r+1}^-=\gamma\cdot{\bf x}_0$ and 
 ${\bf x}_i^-=g_1p_1\cdots g_i\cdot{\bf x}_0, {\bf x}_i^+=g_1p_1\cdots g_ip_i\cdot{\bf x}_0$ for $1\leq i\leq r$. 
We then have:

\begin{lem}\label{horodec}
\noindent Let $\gamma\!=\!g_1p_1\cdots   p_rg_{r+1}$ be the horospherical decomposition of $\gamma$ at height $h$: \linebreak 
\vspace{-4mm}

\noindent (i)  for every $i\in\{ 1,\cdots, r+1\}$  the geodesic segments $[{\bf x}_{i-1}^+, {\bf x}_i^-]$ and $[{\bf x_0}, g_i\cdot{\bf x_0}]$ lie outside the set 
$\bigcup_{g\in\Gamma}g.\mathcal{H}(c)$,  for $c={\rm diam}(\mathcal{K})$;

\noindent (ii) for every $i\in\{ 1,\cdots, r\}$  the geodesic segments $[{\bf x}_0, p_i\cdot {\bf x}_0]$ have length greater than $2(h-c)$ and lie outside the set $\bigcup_{g\in\Gamma}g.\mathcal{H}(h-c)$.

\end{lem}

\noindent{\bf Proof.} 
 Assume $r\geq 1$ and fix $1\leq i\leq r+1$. By construction, each geodesic segment $[{\bf y}_{i-1}^+, {\bf y}_i^-]$ lies outside $\bigcup_{g\in\Gamma}g.\mathcal{H}$.
 Then,  each segment 
$[{\bf x}_{i-1}^+, {\bf x}_i^-]$ lies outside 
$\bigcup_{g\in\Gamma}g.\mathcal{H}(c)$ since $d({\bf x}_{i-1}^+, {\bf y}_{i-1}^+)$ and $d({\bf x}_{i}^-, {\bf y}_{i}^-)$  are both smaller than $c$.   Since $ [{\bf x}_{i-1}^+, {\bf x}_i^-] =g_1p_1\cdots g_{i-1}p_{i-1}\cdot[{\bf x_0}, g_i\cdot{\bf x_0}] $   and the set  $\bigcup_{g\in\Gamma}g.\mathcal{H}(c) $ is $\Gamma$-invariant, the same holds for the segment $[{\bf x_0}, g_i\cdot{\bf x_0}]$.   
%Suppose there exist $i$ and $g$ such that $[{\bf x}_0,g_i\cdot{\bf x}_0]\cap g.\mathcal{H}(a+c)\neq\emptyset$. Translating both sets by $h_{i-1}:=g_1p_1\cdots g_{i-1}p_{i-1}$, the segment $[{\bf x}_{i-1}^+,{\bf x}_i^-]$ intersect should intersect a translate of $\mathcal{H}(a+c)$. Considering the two geodesic segments $[{\bf y}_{i-1}^+,{\bf x}_i^-]$ and $[{\bf x}_{i-1}^+,{\bf x}_i^-]$
%and using convexity property of distance between points on two geodesics (in nonpositive curvature), it appears that $[{\bf y}_{i-1}^+,{\bf x}_i^-]$ should intersect $\mathcal{H}(a+c-diam(\partial\mathcal{C}))$. The same argument with $[{\bf y}_{i-1}^+,{\bf x}_i^-]$ and $[{\bf y}_{i-1}^+,{\bf y}_i^-]$ shows that the latter should intersect $\mathcal{H}(a+c-2diam(\partial\mathcal{C} ))=\mathcal{H}(a)$ which is a contradiction. 
To prove   statement (2), notice that the segment $[{\bf y}_{i}^-, {\bf y}_i^+]$ intersects the set $\bigcup_{g\in\Gamma}g.\mathcal{H}(h)$ and has  endpoints   in $\partial \mathcal{H}$, so that $d({\bf y}_{i}^-, {\bf y}_i^+)\geq 2h$;  one concludes using the facts  that $d({\bf x}_{i}^-, {\bf y}_{i}^-)$ and  $d({\bf x}_{i}^+, {\bf y}_{i}^+)$  are both smaller than $c$ and that
$[{\bf x}_{i}^-, {\bf x}_i^+] =g_1p_1\cdots p_{i-1}g_{i}\cdot[{\bf x_0}, p_i\cdot{\bf x_0}] $.\hfill $\Box$ \pagebreak
\vspace{3mm}
 
 \noindent Moreover,  the distance function is almost additive, with respect to this decomposition : 
\begin{lem} ${}$

\noindent There exists a  constant $C\geq 0$ such that  for all $\gamma\in\Gamma$ with horospherical decomposition $\gamma\!=\!g_1p_1\cdots  g_rp_rg_{r+1}$   at height $h$ :
\vspace{-3mm}

$$d({\bf x}_0,\gamma\cdot{\bf x}_0)\geq \sum_{i=1}^{r+1}d({\bf x}_0,g_i\cdot{\bf x}_0)+\sum_{i=1}^{r}d({\bf x}_0,p_i\cdot{\bf x}_0)- rC.$$
\end{lem}
\noindent{\bf Proof.} 
%The  vertices of the triangle $T({\bf x}_0, g_1\cdot{\bf x}_0, \gamma\cdot{\bf x}_0)$ belong to disjoint $\Gamma$-translates of $\partial\mathcal{H}$ so that the latter has angle at $g_1\cdot{\bf x}_0$ not less than $\pi/2$ and the previous lemma gives a constant $C'$ so that :
%$$\begin{array}{lll}d({\bf x}_0,\gamma\cdot{\bf x}_0)&\geq &d({\bf x}_0,g_1\cdot{\bf x}_0)+d(g_1{\bf x}_0,\gamma\cdot{\bf x}_0)-C'\\
%                                     &=&d({\bf x}_0,g_1\cdot{\bf x}_0)+d({\bf x}_0,p_1g_2\cdots g_{r+1}\cdot{\bf x}_0)-C'.\end{array}$$
%An induction gives the conclusion.
 With the above notations  one gets, for  $C:= 4{\rm diam}(  \mathcal K)$
 \small
 \begin{eqnarray*}
d({\bf x}_0, \gamma \cdot{\bf x}_0)&=& \sum_{i=1}^{r+1} d({\bf y}_{i-1}^+, {\bf y}_{i }^-)+  \sum_{i=1}^r d({\bf y}_i^-, {\bf y}_{i}^+)\\
&\geq & 
 \sum_{i=1}^{r+1} d({\bf x}_{i-1}^+, {\bf x}_{i}^-)+  \sum_{i=1}^r d({\bf x}_i^-, {\bf x}_{i}^+)- 4r {\rm diam}( \mathcal K)\\
 &=& 
  \sum_{i=1}^{r+1} d({\bf x}_0, g_i\cdot {\bf x}_0)+  \sum_{i=1}^r d({\bf x}_0, p_i\cdot {\bf x}_0)- rC. \;\;\;\;\;\;\; \Box 
\end{eqnarray*}
\normalsize

\vspace{2mm}
Now assume $a > \max(h, {\rm diam}(\mathcal{K}))$ and let
$$P_a=\{ p\in P\slash d_0({\bf x}_0,p\cdot{\bf x}_0)\geq 2a\}$$
$$\Gamma_a=\{\gamma\in \Gamma \slash[{\bf x}_0,\gamma\cdot{\bf x}_0]\cap g.\mathcal{H}(a) =\emptyset \mbox{ for all } {g\in\Gamma} \}$$
%  We will need the following   basic result of geometry of pinched Hadamard manifolds and a proof is given in \cite{DPPS1}.
%  \begin{lem} Let $X$ be a Hadamard manifold whose sectional curvature $K$ satisfies $-B^2\leq K\leq -A^2<0$. Then, for every $\alpha_0\in ]0,\pi]$ there exists a constant $C=C(A,B,\alpha_0)$ such that for every triangle whose vertices are ${\bf x},{\bf y},{\bf z}$ and angle $\alpha\geq \alpha_0$ at ${\bf x}$, we have :
%$$d({\bf y},{\bf z})\geq d({\bf y},{\bf x})+d({\bf x},{\bf z})-C.$$
%\end{lem}
 for $t\geq 0$. Let $\gamma\in\Gamma$ with horospherical decomposition $\gamma\!=\!g_1p_1\cdots  g_rp_rg_{r+1}$   at height $h$. By Lemma \ref{horodec}, the geodesic segments $[ {\bf x}_0, g_i\cdot {\bf x}_0],   1\leq i\leq r+1, $ stay outside the perturbed set $\bigcup_{g \in \Gamma}g\cdot \mathcal{ H} (a)$, so that $g_i \in \Gamma_a$ and     $d({\bf x}_0,g_i\cdot{\bf x}_0)=d_0({\bf x}_0,g_i\cdot{\bf x}_0)$.   
 Consequently, the Poincar\'e series of $\Gamma$ for the perturbed metric is 
 \small
\begin{eqnarray*}P_{\Gamma}({\bf x}_0, \delta) \!\!\!\! 
%&=&\ds 1+\sum_{\gamma\neq Id}e^{-\delta d({\bf x}_0,\gamma\cdot{\bf x}_0)}\\
   &\leq  &\! \! \!  \ds 1+ \! \! \! \! \! \! \!   \sum_{\Gamma_{a}\setminus\{ Id\}} \! \! \! \! \! \! \!  e^{-\delta d({\bf x}_0,\gamma\cdot{\bf x}_0)}+
\sum_{p\in P_{a}}    e^{-\delta d({\bf x}_0,p\cdot{\bf x}_0)}
+\sum_{ r\geq 1} \; \sum_{\stackrel{g_1,\cdots, g_{r+1}\in\Gamma_{a}}{_{p_1,\cdots,p_r\in P_a}}}  \! \! \!   \! \! \!  \! \! \! e^{-\delta d({\bf x}_0,g_1p_1\cdots g_rp_rg_{r+1}\cdot{\bf x}_0)}\\
&\leq & \! \! \!   \sum_{\gamma \in \Gamma}e^{-\delta d_0({\bf x}_0,\gamma\cdot{\bf x}_0)}+
\sum_{p\in P_{a}}e^{-\delta d({\bf x}_0,p\cdot{\bf x}_0)}+\sum_{r\geq 1}
\left(e^{C\delta} \! \! \! \sum_{g\in\Gamma_{a}} \! \! \! e^{-\delta d({\bf x}_0,g\cdot{\bf x}_0)}
\sum_{p\in P_{a}}e^{-\delta d({\bf x}_0,p\cdot{\bf x}_0)}\right)^r\\
&\leq &\! \! \!  
   \sum_{\gamma \in \Gamma}e^{-\delta d_0({\bf x}_0,\gamma\cdot{\bf x}_0)}+
\sum_{p\in P_{a}}e^{-\delta d({\bf x}_0,p\cdot{\bf x}_0)}+\sum_{r\geq 1}
\left(e^{C\delta}   \sum_{\gamma \in \Gamma}e^{-\delta d_0({\bf x}_0,\gamma\cdot{\bf x}_0)}
\sum_{p\in P_{a}}e^{-\delta d({\bf x}_0,p\cdot{\bf x}_0)}\right)^r
\end{eqnarray*} 
\normalsize
For $\delta= \delta_P$ the term  $\displaystyle  \sum_{\gamma \in \Gamma}e^{-\delta d_0({\bf x}_0,\gamma\cdot{\bf x}_0)}<+\infty $ since  $\delta_P>N-1$. On the other hand, 
 the group $ P$ being convergent with respect to the new metric $g$, we deduce that 
  $\displaystyle  \sum_{p \in P_a}e^{
-\delta_P d({\bf x}_0), p \cdot {\bf x}_0))} \to 0$ as  $a \to +\infty$; we may thus choose $a$ large enough so that 
$$e^{C\delta_P}   \sum_{\gamma \in \Gamma}e^{-\delta_P d_0({\bf x}_0,\gamma\cdot{\bf x}_0)}
\sum_{p\in P_{a}}e^{-\delta_P d({\bf x}_0,p\cdot{\bf x}_0)}<1
$$ which readily implies $P_{\Gamma}({\bf x}_0,\delta_P)<+\infty$, hence $\delta_P \geq \delta_\Gamma$. 
As  $P$ is a  subgroup  of $\Gamma$, this implies that   $\delta_P = \delta_\Gamma$, hence  $\Gamma$ is a convergent group.\hfill $\Box$

\section{Critical gap property versus divergence}

In this section we start constructing  a  hyperbolic lattice $\Gamma$ of $ \H^2$ which is generated by suitable parabolic isometries,  so that the resulting surface  $\Gamma \backslash \H^2$  has finite volume. \linebreak
In  the Poincar\'e model for the hyperbolic plane, we  choose $r \geq 2$ and $2r$ boundary points $\xi_0 \! =\! \eta_0, \xi_1,\eta_1,\cdots,\xi_r,\eta_r \!=\!\xi_0$ of $\mathbb{S}^1 = \partial \mathbb B^2$ in cyclic order, and consider (uniquely determined) parabolic isometries $p_1,\cdots, p_r$ such that for $1\leq i\leq r$ we have $p_i\cdot\xi_i=\xi_i$ and $p_i\cdot\eta_{i-1}=\eta_i$.  
We remark that  all $\eta_i$ belong to the  $\Gamma$ orbit of the  point $\eta_0$, which is a parabolic fixed point of the isometry  $p_0:=p_r p_{r-1} \cdots  p_1$. 
%$^($\footnote{ Precisely,  for any  fixed $\xi\in\mathbb{S}^1$ there exists a one parameter family of parabolic isometries of the hyperbolic ball $ \mathbb B^2 $ whose fixed point is $\xi$; consequently, for every  points $\eta, \eta' \in \mathbb{S}^1\setminus\{ \xi\}$   there exists a unique parabolic isometry $p$ with fixes $\xi$ and such that $p.\eta=\eta'$. That  makes the choice of the $p_1, \cdots, p_r$ coherent. }$^)$

\begin{property}\label{r+1cups}
The group $\Gamma= \langle p_1, \cdots, p_r\rangle$   is   a free non abelian group over $p_1, \cdots, p_r$. 
The  quotient $\Gamma \backslash \mathbb B^2$  is a   finite surface with $r+1$ cuspidal ends, with a cusp $\mathcal C_i$ for each  parabolic subgroup   $\mathcal P_i=\langle p_i \rangle$ for $i=1,..,r$,  and another cusp   $\mathcal C_{0}$ corresponding  to  the parabolic subgroup $\mathcal P_{0}=\langle p_0 \rangle$  fixing $\xi_0$.
\end{property}

Each element $\gamma \in \Gamma \setminus \{{\rm Id}\}$ can be written in a unique way  as a word  with letters in the alphabet $\mathcal A:= \{p_1^{\pm 1}, \cdots, p_r^{\pm 1}\}$; namely, one gets
\begin{equation}\label{decompositionmathcalA}
\gamma =p_{j_1}^{\epsilon_1}\cdots p_{j_n}^{\epsilon_n} \end{equation}
with $p_{j_1}^{\epsilon_1}, \cdots, p_{j_n}^{\epsilon _n} \in  \mathcal A$, $n \geq 1$  and with adjacent letters which are not inverse to each   other. Such expression with respect to the natural (but not canonical) choice of the alphabet $\mathcal A$ is called a {\em coding} of elements of $\Gamma$. We will call  $j_1$  is  the {\em first index}   of $\gamma$,   denoted by $i_\gamma$; similarly  $j_n$  the {\em last index} and is denoted by $l_\gamma.$ 

\subsection{A new coding for elements of $\Gamma$}
We  code here the elements of $\Gamma$ by blocks, with some admissibility rules to be precised. This new coding is designed to obtain  a contraction property for an operator  that will be introduced  and studied in the next sections  whose restriction to some suitable space of functions present remarkable spectral properties.
\linebreak 
 We first rewrite the decomposition (\ref{decompositionmathcalA}) as follows 
\begin{equation}\label{decompositionpowers}
\gamma =p_{i_1}^{\ell_1}p_{i_2}^{\ell_2}\cdots p_{i_m}^{\ell_m}
\end{equation}
with $m \geq 1, \ell_1, \cdots, \ell_m \in \mathbb Z^*$ and $i_j\neq i_{j+1}$ for $1\leq j<m$.
When all the $\ell_j, 1\leq j\leq m$, belong to $\{\pm 1\}$, one says that $\gamma$ is a  {\it level 1  word}; the set of such words is denoted by $\mathcal W_1$.
Then, we  select all the $\ell_j, 1\leq j\leq m,$ with  $| \ell_j | \geq 2$
%\footnote{Notice that in fact we should need to take also into account the powers of the parabolic isometry $p_0= p_r\cdots p_1$; the modification of the coding should then be more technical and confusing here. We refer to \cite{DPFourier} for the details.  To avoid this complication, we thus  assume that we keep the hyperbolic metric with curvature $-1$ inside the cusp $\mathcal C_0$ (in other words $\tau_0(t)= e^{-t}$ for all $t\geq 0$); it readily follows that, when the constants $a_1, \cdots, a_r$  which define the choice of $g$ are large enough, all the geodesic paths $[{\bf x}_0, q\cdot {\bf x}_0], q \in \widetilde{\mathcal Q}$ remain in the region of $\mathbb B^2$ where the curvature is $-1$, so that the exponential growth of $\widetilde{\mathcal Q}$  is less than  $1$, which will be useful in the sequel.}$^)$;
and write $\gamma$ as  
\begin{equation}\label{newnewdecomposition}
\gamma=p_{j_0}^{l_0}Q_1 p_{j_1}^{l_1}Q_2\cdots p_{j_{k-1}}^{l_{k-1}}Q_kp_{j_k}^{l_{k}}
\end{equation}
where  :
 
$\bullet \quad k\geq 0$

$\bullet \quad  \vert l_1\vert, \cdots, \vert l_{k-1}\vert \geq 2$,

$\bullet \quad \vert  l_0\vert,\, \vert l_k \vert \neq 1$, 

$\bullet \;$   each $Q_j$ is either empty or a level 1  word,   with  $i_{Q_j}\neq j_{i-1}$ and   $l_{Q_j}\neq j_{i}.$  

\noindent The decomposition by blocks (\ref{newnewdecomposition}) is still unique ;  it only uses  letters from the new alphabet
$$
\mathfrak B:= \widehat{\mathcal P}_1\cup \cdots \cup \widehat{\mathcal P}_r\cup \mathcal W_1
$$
 where $\displaystyle \widehat{\mathcal P}_i:=  \{p_i^n/\vert n\vert \geq 2\} $ for  $1\leq i\leq r$,   possibly with
   $p_{i_0}^{l_0}=1$ or  $p_{i_r}^{l_r}=1$. 
 
\noindent We will call {\em blocks} the letters of this new alphabet, and say that a word   $\beta_1\cdots \beta_m$ in the alphabet $\mathfrak B$  is  {\it admissible} if  the last letter of any block  $\beta_i$ is different from  the first one of $\beta_{i+1}$ for  $1\leq i\leq m-1$. 
So, any $\gamma \in \Gamma\setminus\{{\rm Id}\}$ can be written as a finite, admissible word  $\beta_1\cdots \beta_m$ 
on $\mathfrak B$; the  ordered sequence of the $\beta_i$'s is called the {\it $\mathfrak B$-decomposition} of $\gamma$ and 
the number  $m$ of blocks  is  denoted by $\vert \gamma\vert _\mathfrak B$.  
Finally, we denote by $\Sigma_\mathfrak B$ the set of all  finite admissible words with respect to $\mathfrak B$.

\subsection{A new metric in the cusps}
\label{anewmetricinthecusps}

We consider  a fundamental system of horoballs  $\mathcal H_0, \mathcal H_1, \mathcal H_2,$ $\cdots, \mathcal H_{r} $ centered respectively at  the parabolic points $\xi_0, \xi_1, \xi_2, \cdots ,  \xi_r $  and  such that all the horoballs $\gamma\cdot \mathcal H_i$, for $\gamma \in \Gamma, 0\leq i\leq r$, are disjoint or coincide, as in section  \ref{landscape}. \linebreak
%Projections of those horoballs on the hyperbolic surface provide neighborhoods $({\mathcal C}_i)_{0\leq i\leq r}$ of cuspidal ends of that surface on which we modify the metric individually in the following manner. 
Then, we modify the hyperbolic metric in the cuspidal ends $\bar  {\mathcal C}_i = P_i \backslash \mathcal H_i$ as follows. \linebreak
We choose   positive constants $a_0, a_1,    \cdots, a_{r-1}, a_r , b$ and $\eta$,   functions $\tau_0,  \tau_1, \cdots, \tau_{r-1}$ and $ \tau _r$  as in section \ref{cuspidalgeometry} such that $$\omega_{\tau_r} = \max (\omega_{\tau_0}, \cdots, \omega_{\tau_r})>1$$ 
and   we prescribe the profile  $\tau_{i, a_i}$ for the $i$-th cusp   $\bar {\mathcal C}_i$  for $0\leq i\leq r-1$,  and the profile $\tau_{r, a_r,b}$ on the last. dominant cusp ${\mathcal C}_r ={\mathcal C} $. 
This yields   a new surface  $X =(\mathbb B^2, g_{a_0, \cdots, a_r, b})$,  with quotient  $\bar X = \Gamma \backslash X$ of finite area. Since the metric on $X$ depends, in particular, on the value of the parameter $b$, which will play a crucial role in what follows,  
we shall denote the induced distances   on $X$, $\partial X \cong \mathbb{S}^1 $ and the conformal factor respectively by $d_b$, $D_b$ and $| \cdot |_b$; on the other hand, we shall omit the index $b$ in the Busemann function and in the Gromov product,  to simplify notations.
%Notice that, the profiles $\tau_i$ being fixed,  $\Gamma$ has the same critical exponent $\delta$ with respect to  $g_{a_0, \cdots, a_r, b}$, indipendently from the choice of the $a_i$ and $b$.\linebreak
The dependence of $D_b$ on the parameter $b$ is described by   the following lemma,  whose proof can be found in \cite{Pdcds}:

\begin{lem} 
\label{lemmalpha}
Let $b_0>0$ be fixed. There exists $ c \geq 1$ and $\alpha \in]0,1]$ such that  the family of distances 
$(D_b)_{0\leq b \leq b_0}$, are H\"older equivalent; namely   for  all $b \in [0, b_0]$ we have
$$
{1\over c}D_{0}^{1/\alpha}\leq D_{b}\leq c D_{0}^{\alpha}.$$
\end{lem}

 \subsection{Ping-Pong by blocks.}\label{pingpong}
 For any $1\leq i \leq r$, we consider the  sub-arcs $I_i:= [\eta_{i-1}, \eta_i]$ and $I'_i:=  [p_i^{-1}\cdot \eta_{i-1},p_i\cdot  \eta_i]$  of $\mathbb{S}^1 $ containing $\xi_i$. There  exists a
 ping-pong dynamic  between these intervals : namely, for any block $\beta \in \mathfrak B$, we have 
\vspace{1mm}

$\bullet \quad$ if $\beta \in \mathcal W_1$, then 
$\beta \cdot  I'_i \subset  I_{i_\beta}$ for any $i \neq l_\beta$

$\bullet \quad$ if $\beta \in \widehat{\mathcal P}_l$ with $l=i_\beta= l_\beta$, then $\beta \cdot I_i\subset I'_{l }$ for any  $i \neq l  $.
 
 \vspace{1mm} 
\noindent  Moreover, for any $\gamma$ with $\mathfrak B$-decomposition $\beta_1\cdots \beta_m$, we define a compact subset  $K_\gamma \subset \mathbb{S}^1$ as follows:

 \vspace{1mm} 
$\bullet \quad$ $K_\gamma=\cup_{i\neq l_\gamma}I'_i$, if $\beta_m \in \mathcal W_1$

$\bullet \quad$ $K_\gamma=\cup_{i\neq l_\gamma}I_i$, if $\beta_m \in \widehat{\mathcal P}_l$ with $l=l_\gamma$.

  \vspace{1mm} 
\noindent Then, using the fact that  the closure of the sets $I'_i$ and $\partial X\setminus I_i$ are disjoint, one gets:

\begin{lem} \label{Busemann=distance} There exists a constant $C=C(A, \eta)>0$ such that
 $$d_b ({\bf x}_0, \gamma\cdot{\bf x}_0)-C\leq \mathcal B_x (\gamma^{-1}\cdot{\bf x}_0, {\bf x}_0)\leq d_b ({\bf x}_0, \gamma \cdot  {\bf x}_0)$$
 for 
any $\gamma\in \Gamma$ and any  $x\in 
K_\gamma$.
\end{lem}

\noindent 
This lemma implies in particular   the   following  contraction
property (see Prop.2.2 in \cite{BaP2} for a detailed proof):
\begin{lem}\label{contraction1} There exist a real number $r\in ]0, 1[$ 
and a constant 
$C>0$ such that for any $\gamma \in \Gamma$ with length $\vert \gamma\vert _\mathfrak B=k$, one gets 
$$
\forall x \in K_\gamma\qquad\vert \gamma'(x)\vert_0 \leq Cr^k.
$$
\end{lem}
%{\color{red} \`a ecrire juste  pour $D_0$? o\`u pour $b \in [0, b_0]$?}

 \subsection{Coding for the limit points.}\label{coding}

 \noindent An infinite word on the alphabet  ${\mathfrak B}$, i.e. an infinite sequence   $\bar \beta= (\beta_n)_{n\geq 1}$ of elements of $\mathfrak B$ is called {\em admissible} if  any finite subword  $\beta_1\cdots\beta_k$   is admissible ; the set of such words is denoted by  $\Sigma_{\mathfrak B}^+$. 
The  contraction property \ref{contraction1} implies the following fundamental fact:

 \begin{lem} \label{coding}  
 For any $\bar \beta \in \Sigma_{\mathfrak B}^+$,  the sequence $(\beta_1\cdots \beta_n\cdot {\bf x}_0)_{n \geq 0} $ converges  to some point $\pi(\bar \beta)\in \partial X$; the map $\pi: \Sigma_{\mathfrak B}^+\to \partial X$ is one-to-one, 
 and its image  $\pi(\Sigma_{\mathfrak B}^+)$ coincides with the subset $\partial_0 X := 
 \partial X \setminus \bigcup_{i=0}^r \Gamma.\xi_i$.  
 
 \end{lem}

 %def of  $\Sigma_{\mathfrak B}^+ (\gamma)$ = mots qui s'enchainent bien avec $\gamma$

\noindent Notice that,   if $\gamma$ has   ${\mathfrak B}$-decomposition $b_1 \cdots b_m$, then  the subset $K_\gamma$ defined in \ref{pingpong} is the  closure of the subset corresponding, via the coding map $\pi$, to the   the infinite sequences $\bar \beta= (\beta_n)_{n\geq 1}$ such that the {\em concatenation}  $\gamma \ast \bar \beta = (b_1 \cdots b_m   \beta_1 \cdots \beta_i \cdots)$   is admissible. \linebreak
We also set  $J_\gamma:= Cl \{\,\pi(\gamma \ast \bar \beta) \; | \; \bar \beta \in  K_\gamma \} $, that is the closure of  points corresponding to admissible sequences obtained by concatenation with the  ${\mathfrak B}$-decomposition of $\gamma$.

\vspace{1mm}
As indicated previously, this coding by blocks is of interest since the classical shift operator on $\Sigma_{\mathfrak B}^+$ induces  locally, exponentially expanding maps $T^n$ on $\partial_0 X$; the map $T$, described for instance in details in \cite{DPFourier},  has countably many inverse branches, each of them acting by contraction on some subset  of $\partial X$. Namely,  we consider on $\Sigma^+_\mathfrak B$ the  natural shift $\theta$ defined by
$$
 \theta(\bar  \beta):= ( \beta_{k+1})_{k\geq 1}, \quad \forall \, \bar \beta=( \beta_k)_{k\geq 1} \in \Sigma_{\mathfrak B}^+ 
$$
This map induces a transformation $T: \partial_0 X\to \partial_0 X$ via the coding  $\pi$; moreover, 
$T$ can be extended to the whole $\partial X$ by setting, for any $\gamma$ with $\mathfrak B$-decomposition   $\gamma = \gamma_1\gamma_2 \cdots \gamma_n$ 
 \vspace{-2mm}
 
 $$T (\gamma .\xi_i) := \gamma_2 \cdots \gamma_n .\xi_i$$

\noindent and $T(\xi_i):= \xi_i$ for $i=0,...,n$. Then,  for every block $\beta \in \mathfrak B$, the   restriction of $T$  to  $J_\beta$ is the action by $\beta^{-1}$; by the dynamic described above, the inverse branches of the map $T$ have the following property

%\begin{property}\label{contraction branchesinverses}
%There exist $0<r<1$ and a constant $C>0$ such that, for any  sequence $\bar  \beta=( \beta_n)_n \in \Sigma_\mathfrak B^+$, any $k\geq 1$ and  $x, y \in K_{ \beta_k}$  one gets 
%$$D_b( \beta_1\cdots  \beta_k \cdot x ,  \beta_1\cdots  \beta_k\cdot y)\leq C\ r^k D_b(x, y).$$
%\end{property}

\begin{property}\label{contraction branchesinverses}
There exist $0<r<1$ and a constant $C>0$ such that, for any  $\gamma \in  \Gamma$ with   $\vert \gamma \vert_{\mathfrak B} = k$  and for    $x, y \in K_{\gamma}$ we have  
$$
D_0(\gamma \cdot x ,  \gamma \cdot y)\leq C r^k D_0(x, y).
$$
\end{property}
%{\color{red} for any $b \in [0,b_0] $ ?}

\noindent Correspondingly,   for   $\vert \gamma \vert_{\mathfrak B} = k$ and $x, y \in J_\gamma$,  we have the following expanding property:  $$D_0 (T^k x, T^ky) \geq \frac{1}{Cr^k}  D_0(x, y)$$
This property is crucial for the investigation of the spectral properties of the transfer operator, which will be introduced in the next section.

\section{Existence of divergent exotic lattices}
This section is devoted to prove Theorem  \ref{exotic lattice divergent}.  An example of an exotic divergent discrete group has already been constructed in \cite{Pdcds} : it is a Schottky group $\Gamma$ generated  by both  parabolic and hyperbolic isometries, with  $L\Gamma \neq S^1$ and,  consequently, the quotient manifold $\bar X=\Gamma \backslash X$ has infinite volume. 
The strategy in \cite{Pdcds} to obtain both exoticity and divergence is the following : one starts from a free, {\em convergent  Schottky group}  satisfying a ping-pong dynamics, obtained by perturbation of a hyperbolic metric in one cusp; then, it is proved that  the  group becomes divergent  if the   perturbation of the metric is pushed far away in the cusp; finally, it is shown that the divergent group can be made exotic (without losing the divergence property) by pushing  the perturbation of the metric at a suitable height,  by a continuity argument  which is consequence of a  careful  description of the spectrum of  a    transfer operator naturally associated to $\Gamma$.   
Here, we adapt  this approach to obtain a discrete group   with finite covolume. 
\vspace{1mm}

  We start   from the surface $\bar X=\Gamma \backslash X$  with $r+1$ cusps  described in  \ref{anewmetricinthecusps}, with a dominant cusp $\bar {\mathcal C}_r= P_r \backslash \mathcal H_r$  and make $\Gamma$  convergent by  choosing  $a_0,...,a_r \gg0$,  as in Theorem \ref{lattice convergent}. Besides a different coding by blocks (due to the generators which are all parabolic) which gives a slightly different expression for the transfer operator associated to $\Gamma$, the main difficulty here is to show that $\Gamma$  can  be made divergent. This cannot be achieved now by simply  pushing the perturbation  far away  in the cusp, since in our case  $\delta_{P_r}$ is strictly greater than  the critical exponent  of the subset of elements staying in the compact, non-perturbed part;  to obtain the divergence we rather modify $\bar {\mathcal C}_r$ with a profile $\tau_{r, a_r, b} $  which equals the profile  of a cusp with constant curvature metric  $-\omega_{\tau_r}^2$ on a sufficiently large band  of width $b \gg 0$. The, we use the spectral properties of the tranfer operator  to show that there exists a suitable value $b=b^{\ast}$  for which  the lattice becomes simultaneuosly exotic and divergent.
%To ensure at the same time this finiteness property, the generation of the group by parabolic isometries and the freeness property of the group, we need to consider surfaces. We propose a simple coding of the elements of the group inspired from the decomposition  given previously.  

\subsection{On the spectrum of transfer operators}

\noindent  From the analytic point of view, the action of the geodesic flow   is encoded by a {\em transfer operator} associated to the transformation $T$ described above. 
Namely, for any  Borel bounded functions $g\,:\, \partial X\to \mathbb R$ we define

$$
\mathcal L_z\varphi (x)=\sum_{T y=x} e^{-z \mathfrak r(y)}\varphi(y), \;\;\mbox{ } x \in \partial X
$$

\noindent  for a function $ \mathfrak r$ to be defined, depending  on the metric $g_{a_0, \cdots, a_r, b}$,   in particular on the width  $b$ of   the band in the cusp  ${\mathcal C}_r$ where the curvature is $-\omega_{\tau}^2$. We will need to understand precisely the dependance on $b$ when $b\in[0, b_0]$,  so we  stress the dependance on $b$  writing $\mathcal L_{b,z} $ and $ \mathfrak r_b$.\\
\noindent The function $\mathfrak r_b$ we  consider  is given by 
$$\mathfrak r_b: y\mapsto \mathcal B_{\beta^{-1}\cdot y}(\beta^{-1}\cdot{\bf x}_0, {\bf x}_0)=  \mathcal B _{x}(\beta^{-1}\cdot{\bf x}_0, {\bf x}_0) $$  
where $ \mathcal B_x$ is the Busemann function with respect to the metric  $g=g_{a_0, \cdots, a_r, b}$. The idea behind the choice of $\mathfrak r_b$ is to relate the norm  $\vert \mathcal L_{b, s}^{m}1\vert_\infty$ to the Poincar\'e series of $\Gamma$;   on the other hand, the Busemann functions are good approximations of the distance functions $d_b$ by \ref{Busemann=distance},  and are more adapted to express  the iterates of $ \mathcal L_{b, s}$, by the cocycle property.

\noindent To explicitly express the set $T^{-1}(x)$ we remark that,  the alphabet $\mathfrak B$ being countable,  the  pre-images of $x \in\partial X$ by $T$  
%, or equivalently of  $\bar x\in \Sigma_{\mathfrak B}^+$ by $\theta$
are the points $y=\beta \cdot x$,   for those blocks $\beta \in \mathfrak B$ 
 such that $x$ belongs to $K_\beta$, that is :

$\bullet \quad$ if $x \in I'_i$ then $\beta \in \cup_{j\neq i}\widehat{\mathcal P}_j$ or $\beta \in\mathcal W_1$ with  $l_\beta  \neq i$,

$\bullet \quad$ if $x \in I_i\setminus I'_i$ then $\beta \in \cup_{j\neq i}\widehat{\mathcal P}_j$.

\noindent This motivates the introduction  of   {\it weight functions}  $w_{b,z}(\gamma, \cdot): \partial X \rightarrow \mathbb C$, which are defined  for $z\in \mathbb C$ and $ \gamma \in \Gamma $ as
$$
w_{b,z}(\gamma, x)= 1_{K_\gamma}(x) e^{-z \mathfrak r_b (\gamma, x)}
$$

We can now   make explicit the definition of   transfer operators : for  
$b\in [0, b_0], z\in \mathbb C$ and any  bounded Borel function $\varphi: \partial X \to \mathbb C$ we set
\begin{equation}\label{ruelle}
(\mathcal L_{b, z}\varphi) (x)
=\sum_{\beta \in \mathfrak B}  w_{b,z}(\beta, x)\varphi(\beta \cdot x)
 =\sum_{\beta \in \mathfrak B}
1_{K_\beta}(x) e^{-z\mathcal B_x(\beta^{-1}\cdot{\bf x}_0, {\bf x}_0)} \varphi(\beta \cdot x), \;\;\; \forall x \in \partial X  
\end{equation}

\noindent The  operator  $\mathcal L_{b, z}$ is well defined when ${\rm Re} (z)\geq \delta$  and acts on the space $C(\partial X$)  of $\mathbb C$-valued continuous functions on $\partial X$ endowed with the norm $\vert \cdot\vert_\infty$ of the uniform convergence; however, to obtain a quasi-compact operator with good spectral properties, we will consider its restriction to a  subspace  $\mathbb L_\alpha\subset {\rm C}(\partial X)$ of H\"older continuous functions with respect to $D_0$, for $\alpha$ given by Lemma \ref{lemmalpha}.  Namely we let 
$$
\mathbb L_\alpha := \{\varphi\in C (\partial X): \Vert \varphi \Vert = \vert \varphi\vert_\infty+[\varphi]_\alpha <+\infty \}
$$
where
$\displaystyle
[\varphi]_\alpha:= \sup_{\stackrel{x, y \in \partial X}{x\neq y}} {\vert \varphi(x)-\varphi(y)\vert \over D^\alpha_0(x, y)}$; then,    $\mathcal L_{b, z}$ acts on  $\mathbb L_\alpha$ because of the following Lemma 
(which is deduced by Lemma  \ref{lemmalpha} as Lemma III.3 in \cite{BaP2}): 
 \begin{lem}\label{weight}
 Each weight $w_{b,z}(\gamma,\cdot)$ belongs to $\mathbb L_\alpha$ and for any $z \in \mathbb C$, there exists $C=C(z)>0$ such that for any $\gamma \in \Gamma$
 $$
 \Vert w_{b,z} (\gamma, \cdot)\Vert \leq C e^{-{\rm Re}(z) d_b ({\bf x}_0, \gamma\cdot {\bf x}_0)}.
 $$
 \end{lem}

Observe now that   the $w_z(\gamma, \cdot)$ satisfy the following  cocycle relation %\cite{BaP1}
: if   the $\mathfrak B$-decomposition of $\gamma=\gamma_1\gamma_2$ is given by the simple concatenation $\gamma_1\ast \gamma_2$ of the $\gamma_i$, i.e. $\vert \gamma_1\gamma_2\vert _\mathfrak B= \vert \gamma_1 \vert _\mathfrak B+\vert  \gamma_2\vert _\mathfrak B$, then 
 $$w_{b,z} (\gamma_1\gamma_2, x)= w_{b,z} (\gamma_1, \gamma_2  \cdot x)  \cdot  w_{b,z} ( \gamma_2, x).$$ 
This equality leads to the following simple expression of  
the iterates of the transfer operators: for any $n\geq 1$ and $x \in \partial X$
\begin{equation}\label{iterates}
(\mathcal L_{b, z} ^k \varphi )(x)=\sum_{\stackrel{\gamma \in \Gamma}{\vert \gamma\vert _\mathfrak B=n}} 
 w_{b,z}(\gamma, x)\varphi(\gamma \cdot x) =
 \sum_{\stackrel{\gamma \in \Gamma}{\vert \gamma\vert_{\mathfrak B}=k}} 
1_{K_\gamma}(x) e^{-z\mathcal B_x(\gamma^{-1}\cdot{\bf x}_0, {\bf x}_0)} \varphi(\gamma \cdot x) .
\end{equation}

Following  \cite{Pdcds} and using the contraction property 
 in Lemma \ref{contraction1}, we first describe  the spectrum of the operators $\mathcal L_{b, z}$ for real values of the parameter $ z$ greater than or equal to  $\delta :=\omega_\tau /2$:   
\begin{prop}
For any $b \geq 0$ and $s \geq \delta=\omega_\tau/2$, the operator $\mathcal L_{b, s}$ acts both on  $(C(\partial X), \vert\cdot\vert_\infty)$ and  $(\mathbb L_\alpha, \Vert \cdot \Vert)$, with respective spectral radius $\rho_{\infty}(\mathcal L_{b, s})$  and $\rho_\alpha(\mathcal L_{b, s})$.  Furthermore, the operator $\mathcal L_{b, s}$ is quasi-compact\footnote{In other words  its  essential spectral radius on  this space is less than $\rho_\alpha(\mathcal L_{b, s})$} on  $\mathbb L_\alpha$, and:
\begin{enumerate}
\item  $\bullet \quad \rho_\alpha(\mathcal L_{b, s})=   \rho_\infty(\mathcal L_{b, s}),$
\\
$\bullet \quad \rho_\alpha(\mathcal L_{b, s})$ is a simple, isolated eigenvalue  of  $\mathcal L_{b, s}$,
% and  this  eigenvalue  is isolated in the spectrum of  $\mathcal L_{b, s}$ 
\\
$\bullet \quad$the eigenfunction $h_{b, s}$ associated with $ \rho_\alpha(\mathcal L_{b, s})$  is non negative on $\partial X$,

\item for any $s \geq \delta$, the map $b\mapsto  \mathcal L_{b, s}$ is continuous from $\mathbb R^+$ to the space of continuous linear operators on $\mathbb L_\alpha$,

\item   the  function $s \mapsto \rho_\infty(\mathcal L_{b, s})$ is  decreasing on $[\delta, +\infty[$.
\end{enumerate}
\end{prop}

\begin{remark}\label{normalization}
{\em 
Let  $\rho_{b, s} =\rho_\infty(\mathcal L_{b, s})$. Then,   $\mathcal L_{b, s}h_{b, s}=\rho_{b, s}h_{b, s}$, the function $h_{b, s}$ being unique up to a multiplicative constant. By duality, there also exists a unique probability measure $\sigma_{b, s}$ on $\partial X$ such that $\sigma_{b, s} \mathcal L_{b, s}=\rho_{b, s}\sigma_{b, s}$;  the function $h_{b, s}$ becomes uniquely determined   imposing  the  condition  $\sigma_{b, s}(h_{b, s})=1$, which we will assume  from now on.  
}
\end{remark}

\subsection{From convergence to divergence: proof of Theorem \ref{exotic lattice divergent}\label{subsection proof}}
${}$

\noindent Combining  expression (\ref{iterates}) with  Lemma  \ref{Busemann=distance}, one   gets for any $s, b>0$  and $k \geq 0$ 

\begin{equation}
\label{eqnorm}
  \vert \mathcal L_{b, s}^{k}1\vert_\infty  \asymp 
\sum_{\stackrel{\gamma \in \Gamma}{\vert \gamma\vert_{\mathfrak B} =k}} \exp(-s d_b({\bf x}_0,  \gamma\cdot{\bf x}_0)).
\end{equation}

Consequently,  the Poincar\'e series $P_\Gamma(s)$ of $\Gamma$ relatively to $d_b$ and the series 
$\displaystyle \sum_{k\geq 0}\vert \mathcal L_{b, s}^{k}1\vert_\infty $
converge or diverge simultaneously. Following \cite{Pdcds}, we  see  that the function 
$s\mapsto  \rho_\infty( \mathcal L_{b, s})$  is strictly decreasing on $ \mathbb R^+$;  the Poincar\'e exponent  of $\Gamma$ relatively to $d_b$ is then  equal to 
$$
\delta_\Gamma=\sup\Bigl\{ s\geq 0: \rho_\infty( \mathcal L_{b, s})\geq 1\Bigr\}= \inf\Bigl\{ s\geq 0: \rho_\infty( \mathcal L_{b, s})\leq 1\Bigr\}.
$$

This expression of the Poincar\'e exponent in terms of the spectral radius of the transfer operator will be useful to prove  Theorem \ref{exotic lattice divergent}. We first get the

\begin{lem}  \label{lemmab0}
Assume that the profiles $\tau_0, \cdots, \tau_1, \cdots, \tau _r=\tau$  are convergent and satisfy the condition $\omega_\tau = \max (\omega_{\tau_0}, \cdots, \omega_{\tau_r})>2$.
Then there  exist  non negative reals $a_0, \cdots, a_{r}$ and $b_0>0$ such  that 

$\bullet \quad$ the group  $\Gamma$ has exponent ${ \omega_\tau\over 2}$ and is convergent with respect to $ g_{a_0, \cdots, a_r, 0}$;

$\bullet \quad $ the group  $\Gamma$ has exponent $>{\omega_\tau\over 2}$ and is divergent with respect to $g_{a_0, \cdots, a_r, b_0}$.

\end{lem}
\noindent {\bf Proof.  }By the previous section and the choice of the profiles $\tau_i, 0\leq i\leq r$,  we may fix the constants  $a_0, \cdots, a_r$ large enough, in order that the group $\Gamma$ acting on $(X, g_{a_0, \cdots, a_r, 0})$ is a convergent lattice with exponent $\delta={\omega_\tau\over 2} $. Considering  just the contribution of words containing powers of $p_r$, by  the triangle inequality we get  
\vspace{-3mm}

\begin{eqnarray*}
\sum_{\gamma \in \Gamma} e^{-\delta d_b({\bf x}_0, \gamma\cdot {\bf x}_0)}
&\geq& \sum_{k\geq 1}\sum_{\stackrel{Q_1, \cdots, Q_k \in \mathcal W_1 }{\vert l_1\vert , \cdots, \vert l_m\vert \geq 2}}
e^{-\delta d_b({\bf x}_0,  p_r^{l_1}Q_1\cdots p_r^{l_k}Q_k\cdot {\bf x}_0)}\\
&\geq &\sum_{k\geq 1}\left(\sum_{ \vert l \vert \geq 2}
e^{-\delta d_b({\bf x}_0,  p_r^{l}\cdot {\bf x}_0)}\sum_{ Q  \in \mathcal W_1 }
e^{-\delta d_b({\bf x}_0,  Q_ \cdot {\bf x}_0)}\right)^k.
\end{eqnarray*}
When $\vert l\vert $ is large enough,  say $\vert l\vert \geq n_{a_r}\geq 1$, the geodesic segments $[{\bf x}_0, p^l\cdot {\bf x}_0]$ intersect the ``horospherical band'' 
$\mathcal H_{\xi_r} (a_r, a_r+b) := \mathcal  H_{\xi_r} (a_r ) \setminus \mathcal  H_{\xi_r} (a_r+b) $  in which the curvature is  $-\omega_\tau^2$; more precisely, for any $b>0$, there exists $n_b\geq n_{a_r}$ such that:   for any $n_{a_r}\leq l\leq  n_b$,  two points such that there exist integers  $1\leq n_{a_r} < n_b$

When $\vert l\vert $ is large enough,  say $\vert l\vert \geq n_{a_r}\geq 2$, the geodesic segments $[{\bf x}_0, p^l\cdot {\bf x}_0]$ intersect the horoball 
$ \mathcal  H_{\xi_r} (a_r) $; furthermore, for any $b>0$, there exists $n_b\geq n_{a_r}$, with $n_b \to +\infty$ as $b \to +\infty$, such that $[{\bf x}_0, p^l\cdot {\bf x}_0]\cap  \mathcal  H_{\xi_r} (a_r+b) =\emptyset$ when $n_{a_r}\leq \vert l\vert \leq b$ so that the geodesic segment $[{\bf x}_0, p^l\cdot {\bf x}_0]\cap  \mathcal  H_{\xi_r} (a_r) $  remains in a subset of $X$ where the curvature is  $-\omega_\tau^2$. Now, for such $l$, the lenth of the segment $[{\bf x_0}, p^l\cdot {\bf x_0}]$, evaluated either with respect to the metric $g_b$ or to the hyperbolic metric of curvature $-\omega^2_\tau, $ equals  the length of   $[{\bf x}_0, p^l\cdot {\bf x}_0]\cap  \mathcal  H_{\xi_r} (a_r) $ up to $2a_r+c$, for some $c>0$  depending only on the bounds of the curvature; this yields, for $n_{a_r} \leq \vert l\vert \leq n_b$
$$
d_b({\bf x}_0, p^l\cdot {\bf x}_0)\geq   {d_0({\bf x}_0, p^l\cdot {\bf x}_0)\over \omega_\tau}-2(2a_r+c)\geq  { 2\over \omega_\tau}\ln \vert l\vert -2(2a_r+c)$$  
so that,  for $\delta= {\omega_\tau\over 2}$:
$$
\sum_{ \vert l\vert  \geq  2}
e^{-\delta d_b({\bf x}_0, p^{l}\cdot {\bf x}_0)}\geq 
\sum_{ \vert l\vert =n_{a_r}}^{n_b}
e^{-\delta d_b({\bf x}_0, p^{l}\cdot {\bf x}_0)}\succeq \sum_{ \vert l\vert =n_{a_r}}^{n_b}\left({1\over \vert l\vert ^{2/\omega_\tau}}\right)^{\omega_\tau/2} \to +\infty \quad {\rm as} \quad b\to +\infty.
$$
So,  there exists $b_0\geq 0$ such that   $\displaystyle \left(\sum_{ \vert l\vert  \geq 2}
e^{-\delta d_b({\bf x}_0, p^{l}\cdot {\bf x}_0)}\sum_{ Q  \in \mathcal W_1 }
e^{-\delta d_b({\bf x}_0,  Q \cdot {\bf x}_0)}\right)>1$  as soon as  $b \geq b_0$; in particular, by monotone convergence type argument, one gets 
$$ \left(\sum_{ \vert l\vert  \geq 2}
e^{-s d_{b_0}({\bf x}_0, p^{l}\cdot {\bf x}_0)}\sum_{ Q  \in \mathcal W_1 }
e^{-s d_{b_0}({\bf x}_0,  Q \cdot {\bf x}_0)}\right)>1$$
 for some $s>{\omega_\tau\over 2}$. This  ensures that  $\delta_\Gamma\geq s> {\omega_\tau\over 2}=\delta_{\mathcal P_r} = \max \{ \delta_{\mathcal P_i} \, | \,  0\leq i\leq r \} $,  and  $\Gamma$ is divergent with respect  to the metric $g_{a_0, \cdots, a_r, b}$ by the critical gap  property recalled in the introduction.  $\Box$

\vspace{1cm}

\noindent {\bf Proof of Theorem \ref{exotic lattice divergent}}.  Recall that $\delta=\omega_\tau/2$. 
 Since  $\rho_\alpha( \mathcal L_{b, \delta})$ is an eigenvalue of $\mathcal L_{b, \delta}$ which is isolated in the spectrum of  $\mathcal L_{b, \delta}$, the function $b\mapsto \rho_\alpha(\mathcal L_{b, \delta})$ has the same regularity as $b \mapsto \mathcal L_{b, \delta}$. 
% {\color{red}  plutot $\mathcal L_{b, s}$ partout? sinon il faut dire $\delta=?$ \\
% introduire $\delta_{\Gamma, b}$ pour la suite} \\
For $b_0$ given by Lemma  \ref{lemmab0}, we have   $\rho_\alpha(\mathcal L_{0, \delta})=\rho_\infty(\mathcal L_{0, \delta})\leq 1$
 and  $\rho_\alpha(\mathcal L_{b_0, \delta})=\rho_\infty(\mathcal L_{b_0, \delta})\geq 1$; thus, there  exists $b_*\in [0, b_0]$ such that  $\rho_w(\mathcal L_{b^*, \delta})=  \rho_\infty(\mathcal L_{b^*, \delta})=1$.
 Since the function  $s \mapsto \rho_\infty(\mathcal L_{b^*, s})$ is strictly decreasing on $[\delta, +\infty[$, one gets $\rho_\infty(\mathcal L_{b^*, s})<1$ as soon as  $s>\delta$. For such values of $s$, the Poincar\'e series $P_\Gamma(s)$ of $\Gamma$  relatively to the metric $g_{b^*}$ thus converges,
this implies that  its Poincar\'e exponent $\delta_{\Gamma,  b^*}$ is  $\leq \delta$; we have in fact $\delta_{\Gamma,  b^*}=\delta$ since $\delta_{\langle p_r\rangle}=\delta$ and $p_r \in \Gamma$.
Finally, the eigenfunction $h_{b^*, \delta}$ of $\mathcal L_{b^*, \delta}$  associated with $\rho_\alpha(\mathcal L_{b^*, \delta})$ being non negative on $\partial X$, one gets  $h_{b^*, \delta}\asymp 1$ and so
 $$
 \sum_{k\geq 0} \vert \mathcal L_{b^*, \delta}^{k} 1\vert_\infty  \asymp \sum_{k\geq 0} \vert \mathcal L_{b^*, \delta}^{k} h_{b^*, \delta}\vert_\infty =   \sum_{k\geq 0} \vert  h_{b^*, \delta}\vert_\infty=\infty$$ 
which  implies, by (\ref{eqnorm}),  that $\Gamma$ is divergent with respect to the metric $g_{b^*}$. $\Box$

\section{Counting for some divergent exotic lattices} 

We consider here the group $\Gamma=\langle p_1, \cdots, p_r\rangle, r\geq 2, $ where the $p_i, 1\leq i\leq r$, satisfy the conditions given in Property  \ref{r+1cups}; we endow $X$ with a  metric of the form
$g= g_{a_0, \cdots, a_r, b^*}$ as defined in subsection \ref{anewmetricinthecusps}, 
where 
 profiles $\tau_0, \tau_1, \cdots, \tau _r$  satisfy conditions  (\ref{conditionsur1}), (\ref{conditionsur2}), (\ref{conditionsur3})  and   constants  $a_0, a_1,  a_2, \cdots, a_r$ and $b^*$   are  suitably chosen in such a way

 $\bullet \ {\bf H_0}$ : {\it the group $\Gamma$ is exotic and divergent  with respect to $g$, with Poincar\'e exponent  $\delta= \delta_\Gamma= {\omega\over 2}$ with $\omega=\max (\omega_{\tau_0}, \cdots, \omega_{\tau_r})>2$.}

 $\bullet \ {\bf H_1}$ :  {\it the exponential growth of the set  $\{\beta \in \mathcal W_1\slash d({\bf x}_0, \beta\cdot{\bf x}_0)\leq R\}$ as $R\to +\infty $  is strictly less than $\delta$.}

$\bullet \ {\bf H_2}$ :  {\it there exist   $\kappa \in ]1/2, 1[$ and a slowly varying function $^($\footnote{A function $L(t)$ is said to be ''slowly varying'' or ''of slow growth'' if it is positive, measurable and $L(\lambda t)/L(t)\to  1$  as $t\to +\infty$ for every  $\lambda>0$.}$^)$  $L$    such that 
\begin{equation}\label{parabolic r}
\sum_{p\in \mathcal P_r\slash  d({\bf x}_0, p\cdot {\bf x}_0)>t}e^{-\delta  d({\bf x}_0,p\cdot\bf{o})}\quad 
  {\buildrel{{\scriptscriptstyle t \to +\infty }}\over \sim}\quad 
  {L(t)\over t^\kappa}.
 \end{equation}}

$\bullet\  {\bf H_3 }$ :  {\it the groups $\mathcal P_l, 0\leq l \leq r-1,$ satisfy the condition
\begin{equation}\label{autres paraboliques}
\sum_{p\in \mathcal P_l\slash  d({\bf x}_0, p\cdot{\bf x}_0)>t}e^{-\delta  d({\bf x}_0,\gamma\cdot\bf{o})} \quad = \quad o\Bigl({L(t)\over t^\kappa}\Bigr).
\end{equation}
}

Since $b^*$ is now fixed, we will omit it it the sequel and set $d=d_{b^*}$.

Under those hypotheses, the subgroup $\mathcal P_r$  corresponding to the cusp ${\mathcal  C}_r$ has   a dominant influence on the behavior of the orbital function of $\Gamma$. Hypothesis  {\bf H2}   is inspired by probability theory, it corresponds to a  {\it heavy tail  conditions} satisfied by random walks which have been    intensively investigated  \cite{E}.
Theorem \ref{counting} will be a direct consequence of the following one:

\begin{theo}\label{countingkappa} Let  $\Gamma$ be   a finite area Fuchsian group  satisfying Property \ref{r+1cups} and let us endow $\mathbb X$ with a metric $g$ such that hypotheses ${\bf H_0, H_1, H_2}$ and $ {\bf H_3}$ hold.  Then,  for any $1\leq j\leq r$, any fixed $x_j\notin I_j$, there exists $C_j>0$ such that  
\begin{equation}\label{busemanrenewaltheorem}
\sharp\{\gamma\in \Gamma_j \slash  \mathcal B_{x_j}(\gamma^{-1}\cdot{\bf x}_0, {\bf x}_0)\leq R \}
\quad {\buildrel{{\scriptscriptstyle R \to +\infty }}\over \sim}
\quad   C_j{  e^{\delta_\Gamma R} \over  R^{1-\kappa}L(R)}.
\end{equation}
where $\Gamma_j$ denotes the set of $\gamma \in \Gamma$ with last letter $j$.
\end{theo}
Indeed, since $x_j \notin I_j$ is fixed, there exists $c(x_j)>0$ such that  
$$
d({\bf x}_0, \gamma \cdot {\bf x}_0)-c(x_j)\leq 
\mathcal B_{x_j}(\gamma^{-1}\cdot{\bf x}_0, {\bf x}_0)
\leq  
d({\bf x}_0, \gamma \cdot {\bf x}_0)$$
 for any $\gamma \in \Gamma_j$ and (\ref{busemanrenewaltheorem}) thus implies
$$
\sharp\{\gamma\in \Gamma_j \slash d( {\bf x}_0, \gamma\cdot{\bf x}_0)\leq R \}
\quad {\buildrel{{\scriptscriptstyle R \to +\infty }}\over \asymp}
\quad   {  e^{\delta_\Gamma R} \over  R^{1-\kappa}L(R)}.
$$
Theorem   \ref{counting}  follows summing over  $j\in \{1, \cdots, r\}$.

Let us make some remarks :
\begin{enumerate} 
\item We have seen in the previous section how one may choose profiles $\tau_0, \cdots, \tau_r$ and parameters $a_0, \cdots, a_r$ and $b^*$ in such a way hypothesis ${\bf H_0}$ holds.
\item 
 Hypothesis ${\bf H_1}$ holds when $\omega:=\max (\omega_{\tau_0}, \cdots, \omega_{\tau_r})>2$, $a_0=+\infty$ and the constants $a_1, \cdots, a_n$ are chosen large enough  in such a way  all the geodesic paths $[{\bf x}_0, \beta\cdot{\bf x}_0], \beta \in \mathcal W_1,$ remain  in a closed subset of $X$ where the curvature is $-1$; the exponential growth of $\{\beta \in \mathcal W_1\slash d({\bf x}_0, \beta\cdot{\bf x}_0)\leq R\}$ as $R\to +\infty $  is  less or equal to $1$ and then strictly less than $\delta=\omega/2$.
\item Hypothesis ${\bf H_2}$ holds in particular when 
$$
d({\bf x}_0, p_r^{n}\cdot {\bf x}_0)={2\ln n+ 2(1+\kappa)(\ln \ln n+O(n))\over \omega}
$$
for some $\omega>2$.
 This equality involves only the asymptotic geometry on the cusp $\mathcal C_r$ as it is equivalent to prescribe a profile. Hence it is compatible with any choice of the parameter $b$.
 The critical exponent of $\mathcal P_r$ is thus $\delta={\omega/2}$
and one gets, as $t\to +\infty,$
\begin{eqnarray*}
\sum_{p\in \mathcal P_r\slash  d({\bf x}_0, p\cdot {\bf x}_0)>t}e^{-\delta  d({\bf x}_0,p\cdot\bf{o})}&=&\sum_{n \in \mathbb N\slash  d({\bf x}_0, p_r^n\cdot {\bf x}_0)>t}e^{-\delta  d({\bf x}_0,p_r^n\cdot\bf{o})}\\
&\asymp& \sum_{n>{e^{\omega t/2}\over t^{1+\kappa}}} {1\over n(\ln n)^{1+\kappa}}\\
&\asymp& \int_{e^{\omega t/2}\over t^{1+\kappa}}^{+\infty} {du\over u(\ln u)^{1+\kappa}}\\&\asymp&{1\over t^\kappa}.
\end{eqnarray*}
\item The condition $\kappa \in ]1/2, 1[$ readily implies 
\begin{eqnarray*}
\sum_{p\in  \mathcal P_r} d({\bf x}_0, p\cdot {\bf x}_0)e^{-\delta d({\bf x}_0,p\cdot {\bf x}_0)}
&\asymp&
 \sum_{N\geq 1} N\Bigl( 
\sum_{p\in  \mathcal P_r\slash  N< d({\bf x}_0, p \cdot {\bf x}_0)\leq N+1}e^{-\delta  d({\bf x}_0, p \cdot {\bf x}_0)}
\Bigr)\\
&\asymp&
 \sum_{N\geq 1}\Bigl( 
\sum_{p\in  \mathcal P_r\slash   d({\bf x}_0, p \cdot {\bf x}_0)>N}e^{-\delta  d({\bf x}_0, p\cdot\bf{o})}
\Bigr)\\
&\asymp& 
 \sum_{N\geq 1}{L(N)\over N^\kappa}=+\infty.
\end{eqnarray*}
 By Theorem B in \cite{DOP}, it readily follows that the Bowen-Margulis measure associated with $\Gamma$ is infinite. 
 \item 
 Observe at last that hypothesis ${\bf H_3}$ is satisfied in particular when the  gap property $\delta_{\mathcal P_l}<\delta$ hold for any $0\leq l\leq r-1$.
\end{enumerate}

 We present in this section  the main steps of the proof of Theorem \ref{countingkappa}; we refer to \cite{DP}  for details. 
For any $R\geq 0$, let us denote $W_j(R, \cdot)$ the measure on $\mathbb R$ defined by: for any Borel non negative function $\psi: \mathbb R\to \mathbb R$
 $$
 W_j(R, \psi):= \sum_{\gamma \in \Gamma_j} e^{-\delta \mathcal B_{x_j}(\gamma^{-1}\cdot{\bf x}_0, {\bf x}_0)}\psi(\mathcal B_{x_j}(\gamma^{-1}\cdot{\bf x}_0, {\bf x}_0)-R).
 $$
 One gets $0\leq W_j(R, \psi) <+\infty$ when $\psi$ has a compact support in $ \mathbb R $ since the group $\Gamma$ is discrete, furthermore $  \sum_{j=1}^rW_j(R, \psi)=e^{-\delta R}v_\Gamma(R)$ when $\psi(t) = e^{\delta t} 1_{t\leq 0} $. If one proves that for any non negative  and continuous function $\psi$ with compact support  and such that  $\int_{\mathbb R^+}\psi(x)dx >0$
 \begin{equation} \label{convergence vague}
 W_j(R, \psi) \quad 
{\buildrel{{\scriptscriptstyle R\to  +\infty }}\over \sim}\quad  {C_j\over R^{1-\kappa}L(R)} \int_{\mathbb R}\psi(x) dx, 
\end{equation}
 this convergence will also hold for non negative functions with compact support in $\mathbb R $ and whose discontinuity set has $0$ Lebesgue measure; 
Theorem \ref{countingkappa} follows since $v_\Gamma(R)= e^{\delta R}\sum_{j=1}^r\sum_{n\geq 0}W_j(R, e^{\delta  t}1_{]-(n+1), -n]}(t))$.
From now on, we fix  a continuous function $\psi: \mathbb R\to \mathbb R^+$ with compact support; one gets, for $1\leq j\leq r$ fixed 
$$
 W_j(R, \psi)
=  \sum_{k\geq 0}
\Bigl(
 \sum_{\gamma \in \Gamma_j\slash  \vert \gamma\vert_\mathfrak B =k}e^{-\delta \mathcal B_{x_j}(\gamma^{-1}\cdot{\bf x}_0, {\bf x}_0)}\psi(\mathcal B_{x_j}(\gamma^{-1}\cdot{\bf x}_0, {\bf x}_0)-R) \Bigr).
$$
Notice that for   $\gamma\in \Gamma_j$ with $\mathfrak B$-decomposition $\gamma= \beta_1\cdots \beta_k$, one gets  
\begin{eqnarray*} 
\mathcal B_{x_j}(\gamma^{-1}\cdot{\bf x}_0, {\bf x}_0)&=&
\mathfrak r(\beta_1\cdots \beta_k  \cdot x_j)+
\mathfrak r(\beta_2\cdots \beta_k  \cdot x_j)+\cdots
+\mathfrak r(\beta_k  \cdot x_j) \\
&=&
\mathfrak r(y)+\mathfrak r(T\cdot y)+\cdots+
\mathfrak r(T^{k-1}\cdot y)= S_k{\mathfrak r}(y)
\end{eqnarray*}
where one denotes  $y:= \gamma\cdot x_j$, so that
\begin{equation}\label{Wensomme}
W_j(R, \psi)
=  \sum_{k\geq 0}
\sum_{y\in \partial X\slash T^k\cdot y= x_j} e^{-\delta S_k\mathfrak r(y)}\psi(S_{k}\mathfrak{r}(y)-R).
\end{equation}

 By a classical argument in probability theory due to Stone (see for instance \cite{E}), it suffices to check that the convergence (\ref{convergence vague}) holds when $\psi$  has a $C^\infty$ Fourier transform with compact support: indeed, the test function  $\psi$ may vary  in the set $\mathcal H$ of functions of the form $\psi(x)= e^{itx}\psi_0(x)$ where $\psi_0$ is an integrable and strictly positive function  on $\mathbb R$ whose Fourier transform is $C^\infty$ with compact support.  When $\psi\in \mathcal H$, one can use 
  the inversion Fourier formula and write
$ 
\psi(x) ={1\over 2\pi} \int_{\mathbb R} e^{-itx}\hat{\psi}(t)dt,
$
so that, for any $0<s<1$
\begin{eqnarray}\label{Wfourier}
W_j(s, R, \psi)
&:=& \sum_{k\geq 0} s^k
 \sum_{y\in \partial X\slash T^k\cdot y= x_j} e^{-\delta S_k\mathfrak r(y)}\psi(S_{k}\mathfrak{r}(y)-R)\notag\\
 &=& \sum_{k\geq 0} s^k
 \sum_{y\in \partial X\slash T^k\cdot y= x_j} {1\over 2\pi} \int_{\mathbb R}e^{itR}e^{-(\delta+it) S_k{\mathfrak r}(y)}\hat{\psi}(t)dt \notag\\
 &=&{1\over 2\pi} \int_{\mathbb R}e^{itR}\hat{\psi}(t)\Bigl(\sum_{k\geq 0} s^k\sum_{y\in \partial X\slash T^k\cdot y= x_j} e^{-(\delta+it) S_k{\mathfrak r}(y)}\Bigr)dt
\notag\\
&=&{1\over 2\pi} \int_{\mathbb R}e^{itR}\hat{\psi}(t)\Bigl(\sum_{k\geq 0} s^k \mathcal L_{ \delta+it}^k1(x_j)\Bigr)dt\notag
\\
&=&{1\over 2\pi} \int_{\mathbb R}e^{itR}\hat{\psi}(t)\left(I-s\mathcal L_{\delta+it} \right)^{-1}1(x_j)dt
\end{eqnarray}
where, for any $z \in \mathbb C$, the operator $\mathcal L_{z}$ is the transfer operator associated to the function ${\mathfrak r}$ formally defined in the previous section by: for any Borel bounded function $\varphi: \partial X \to \mathbb C$ and any $x \in \partial X$
$$
\mathcal L_{ z}\varphi(x):=\sum_{y\in \partial X\slash T\cdot y= x}e^{-z{\mathfrak r}(y)}\varphi(y).
$$
We know that the operators $ \mathcal L_{z}, {\rm Re}(z)\geq \delta$, are bounded and quasi-compact on the space $\mathbb L_\omega(\partial X)$ of  H\"older continuous  function on $(\partial X, D_0)$. 
The   subsection  \ref{controlspectrum} is  devoted  to the control of the peripherical spectrum of  $\mathcal L_{\delta+it} $ on $ \mathbb L_\omega$. In subsection \ref{local}, we describe the local expansion of the dominant eigenvalue. Atlast we  achieve the proof using arguments coming from  renewal theory (subsection \ref{renewal}).

\subsection{ \label{controlspectrum} 
 The spectrum of $\mathcal L_{\delta+it} $ on $ \mathbb L_\omega$ }

First  we need to control the spectrum of $\mathcal L_{z}$ when $z=\delta+it, t \in \mathbb R$.
By Lemma \ref{weight}, the operators $\mathcal L_z$ are bounded on $ \mathbb L_\omega$ when ${\rm Re}(z) \geq \delta$; the spectral radius of $\mathcal L_z$ will be denoted $\rho_\omega(z)$ throughout this section. In the following Proposition, we describe its  spectrum   on  $ \mathbb L_\omega$ when ${\rm Re} (z)=\delta$. 
\begin{prop} \label{spectrum}There exist  $\epsilon_0>0$   and $\rho_0\in ]0, 1[$ such that, for any $t \in \mathbb R$ with modulus less than $\epsilon_0$, the spectral radius $\rho_\omega(\delta+it)$ of $\mathcal L_{\delta+it} $ is $ >\rho_0$ and 
 the operator $\mathcal L_{\delta+it} $ has a unique  eigenvalue $\lambda_t$ of modulus  $\rho(\delta+it)$, which is simple and closed to $1$, the rest of the spectrum being included in a disc of radius $\rho_0$.
 
 Furthermore,  for any $A>0$,  there exists $\rho_A\in ]0, 1[$ such that   $\rho_\omega(\delta+it)\leq \rho_A$ for any $t \in \mathbb R$ such that $\epsilon_0\leq \vert t\vert \leq A.
$
\end{prop}

 \begin{notation}  
We denote $\sigma$ the unique probability measure on $\partial X$ such that $\sigma \mathcal L_\delta =\sigma$ and $h$ the element of $\mathbb L_\omega$ such that $\mathcal L_\delta h=h$ and $\sigma(h)=1$. 
\end{notation}

\noindent {\bf Proof of Proposition \ref{spectrum}} This is exactly the same proof that the one presented in \cite{BaP1}(Proposition 2.2) and \cite{DP}:   the operators  $\mathcal L_{\delta+it} $ are quasi-compact on $ \mathbb L_\omega$ and it is sufficient to control their peripherical spectrum. When $t$ is closed to $0$, we use the perturbation theory to conclude that the spectrum of $\mathcal L_{\delta+it} $ is closed to the one of $\mathcal L_{\delta} $: it is thus necessary to prove that  the map $t\mapsto \mathcal L_{\delta+it} $ is continuous on $\mathbb R$. 
The following Lemma is devoted to precise the type of continuity of this function.
\begin{lem} \label{regularitepourperturbationoperateurs} Under the hypotheses ${\bf H_1}, {\bf H_2}$ and ${\bf H_3}$, there exists a constant $C>0$ such that
$$
\Vert \mathcal L_{\delta+it'} - \mathcal L_{\delta+it} \Vert \leq C \vert t'-t\vert^\kappa L\left({1\over \vert t'-t\vert}\right)
$$
\end{lem}
\noindent {\bf Proof. }
We will use the following classical fact (\cite{Fe} p.272 and \cite{E} Lemmas 1 $\&$ 2):
\begin{lem}
Let $\mu$ be a probability measure on $\mathbb R^+$, set $F_\mu(t):=\mu[0, t]$  and $m(t):= \int _0^t(1-F_\mu(s))ds$ and assume that there exist $\kappa \in \mathbb R$ and a slowly varying function $L$ such that
$1-F_\mu(t)\sim  {L(t)\over t^\kappa}$ as $t\to +\infty$. One thus gets
\begin{equation} \label{muvariationlente}
\lim_{t\to +\infty}{t(1-F_\mu(t))\over m(t)}=1-\kappa 
\end{equation}
and the characteristic function $\displaystyle \hat{\mu}(t):= \int_0^{+\infty} e^{itx}\mu(dx)$ of $\mu$ has the following local expansion as $t\to 0^+$
\begin{equation}\label{localexpansion}
\hat{\mu}(t)=1-e^{-i\pi {\kappa\over 2}}\Gamma(1-\kappa) t^\kappa L\Bigl({1\over t}\Bigr) (1+o(t)).
\end{equation}
\end{lem}
\noindent Noticing that $\tilde{m}(t):= \int_0^tx\mu(dx)= m(t)-t(1-F_\mu(t))$ one also gets $\tilde{m}(t)\sim \kappa m(t)$ as  $t\to +\infty$; furthermore, decomposing  $\int_0^{+\infty}\vert e^{itx}-1\vert\ \mu(dx)$ as 
$\int_{[0, 1/t]}\vert e^{itx}-1\vert\ \mu(dx)+\int_{]1/t, +\infty[}\vert e^{itx}-1\vert\ \mu(dx)$
and applying the previous estimations, one gets, for any $t\in \mathbb R$
$$
 \int_0^{+\infty}\vert e^{itx}-1\vert\ \mu(dx) \preceq t^\kappa L(1/t).
$$
We now apply (\ref{muvariationlente}) with  the probability measures $\mu_l, 1\leq l\leq r$, on $\mathbb R^+$  defined by $\mu_l:= c_l \sum_{p \in \mathcal P_l}\delta_{d({\bf x}_0, p\cdot {\bf x}_0)}$ where $c_l>0$ is some normalizing constant.
As a direct consequence, under   hypotheses   ${\bf H_1}, {\bf H_2}$ and ${\bf H_3}$, one gets (up to a modification of $L$ by multiplicative constant)
\begin{equation} \label{tildempourd}\sum_{p\in \mathcal P_r\slash  d({\bf x}_0, p \cdot{\bf x}_0)\leq t} d({\bf x}_0, p  \cdot {\bf x}_0)e^{-\delta d({\bf x}_0, p  \cdot{\bf x}_0)}\quad 
{\buildrel{{\scriptscriptstyle t\to  +\infty }}\over \sim}\quad  {t^{1-\kappa}\over 1-\kappa}L(t).
\end{equation}
and 
\begin{equation} \label{fourierpourdetL}
\sum_{p\in \mathcal P_r} \vert e^{it  d({\bf x}_0, p  \cdot{\bf x}_0)}-1\vert \times e^{-\delta d({\bf x}_0, p \cdot{\bf x}_0)}\preceq t^{\kappa }L(1/t);
\end{equation}
similarly, for $
1\leq l <r$, one has 
\begin{equation} \label{fourierpourdetK+1/K+L}
\sum_{p\in \mathcal P_l}\vert e^{it  d({\bf x}_0, p \cdot{\bf x}_0)}-1\vert \times e^{-\delta d({\bf x}_0, p  \cdot{\bf x}_0)}= t^{\kappa }L(1/t) o(t). 
\end{equation} 
In the same way, by lemma  \ref{Busemann=distance}, for   any   $x \in \partial X\setminus I_r$, one gets 
\begin{equation} \label{tildempourb}
\sum_{
p\in \mathcal P_r\slash   \mathfrak r (p\cdot  x)\leq t
}  e^{-\delta  {\mathfrak r} (p \cdot x))}{\asymp}  {L(t)\over t^{\kappa}},\quad 
\sum_{p\in \mathcal P_r\slash  {\mathfrak r} (p  \cdot x))> t}  
{\mathfrak r} (p  \cdot x)
e^{-\delta {\mathfrak r} (p  \cdot x)} {\asymp}  t^{1-\kappa} L(t),
\end{equation} 
\begin{equation} \label{fourierpourdetL}
\sum_{p\in \mathcal P_r} \vert e^{it  {\mathfrak r} (p \cdot x)}-1\vert \times e^{-\delta {\mathfrak r} (p  \cdot x)} {\preceq} t^{\kappa }L(1/t)
\end{equation}
and, for   $ 1\leq l\leq r-1$ and $x \in \partial X\setminus I_l$ 
\begin{equation} \label{fourierpourbetK+1/K+L}
\sum_{p\in \mathcal P_l} \vert e^{it  {\mathfrak r} (p \cdot x)}-1\vert \times e^{-\delta {\mathfrak r} (p \cdot x)} = \ t^{\kappa }L(1/t)o(t).
\end{equation}

\noindent Similarly, since the exponential growth of the set $\{\beta \in \mathcal W_1\slash d({\bf x}_0, \beta\cdot{\bf x}_0)\leq R\}$ is $<\delta$, we also have, for any $1\leq l\leq r$  and 
$x\in \cup_{l'\neq  l}I'_{l'}$
\begin{equation} \label{fourierpourbetK+1/K+L}
\sum_{\stackrel{\beta \in \mathcal W_1}{l_\beta=l}} \vert e^{it \mathfrak r(\beta    \cdot x)}-1\vert \times e^{-\delta   {\mathfrak r}(\beta\cdot x)}  \ t^{\kappa }L(1/t)o(t).
\end{equation}
Noticing that for any $\beta \in \mathfrak B$ and $x \in K_\beta$ 
$$
\vert w_{\delta+it}(\beta, x)-w_{\delta+it'}(\beta, x)\vert\leq 
\vert   e^{i(t-t')  \cdot x))}-1\vert \times e^{-\delta   {\mathfrak r}(\beta\cdot x)} 
$$
one readily gets, combining the above   estimations  (\ref{tildempourb}) to (\ref{fourierpourbetK+1/K+L}) all together
\begin{equation}\label{normeinfinie}
 \big\vert \mathcal L_{\delta+it} -\mathcal L_{\delta+it'} \big\vert_\infty\leq \sum_{\beta \in \mathfrak B}
\vert w_{\delta+it}(\beta, \cdot)-w_{\delta+it'}(\beta, \cdot)\vert _\infty\  \preceq  \vert t-t'\vert^{\kappa }L\left({1\over \vert t-t'\vert}\right).
\end{equation}
  Now, for any $\beta \in \mathfrak B$ and $x, y \in K_\beta$, one gets 

$\displaystyle 
\Big\vert
 \bigl(w_{\delta+it}(\beta, x)-w_{\delta+it'}(\beta, x)\bigr)
-
 \bigl(w_{\delta+it}(\beta, y)-w_{\delta+it'}(\beta, y)\bigr)\Big\vert\\
$
\begin{eqnarray*}
\quad &\leq &
  e^{-\delta   {\mathfrak r}(\beta\cdot x)}\times \Big\vert   \left( e^{i(t-t') {\mathfrak r}(\beta\cdot x)}-1\right)- \left( e^{i(t-t'){\mathfrak r} (\beta\cdot  y)}-1\right)\Big\vert\\
  & &\quad +\quad \Big\vert   e^{-\delta   {\mathfrak r}(\beta\cdot x)}-e^{-\delta {\mathfrak r} (\beta\cdot y)}\Big\vert \times \Big\vert   e^{i(t-t'){\mathfrak r} (\beta\cdot  y)}-1\Big\vert\\
  &\leq&
 e^{-\delta   {\mathfrak r}(\beta\cdot x)}\times \Big\vert  e^{i(t-t')({\mathfrak r} ( \beta \cdot x)-{\mathfrak r}  (\beta\cdot y))}-1\Big\vert\\
  & &\quad +\ 
  e^{-\delta {\mathfrak r}  (\beta\cdot y)}  \Big\vert  e^{\delta({\mathfrak r}  (\beta\cdot  y)- {\mathfrak r}  (\beta\cdot  x))}-1\Big\vert
 \times \Big\vert   e^{i(t-t'){\mathfrak r}  (\beta\cdot y)}-1\Big\vert 
  \\
  &\preceq& \Bigl(e^{-\delta   {\mathfrak r}(\beta\cdot x)} 
   [r\circ  \beta ] \times   \vert t-t'\vert    +\  e^{-\delta {\mathfrak r}  (\beta\cdot y)} [r \circ  \beta  ] \times 
 \Big\vert   e^{i(t-t'){\mathfrak r}  (\beta\cdot y)}-1\Big\vert  
\Bigr)D_0(x, y)\end{eqnarray*}
so that, as above
$\displaystyle 
\sum_{\beta \in \mathfrak B}[w_{\delta+it}(\beta, \cdot)-w_{\delta+it'}(\beta, \cdot)]
 \preceq 
 \vert t-t'\vert^{\kappa }L\left({1\over \vert t-t'\vert}\right).
$
  We achieve  the proof of the Lemma  combining this last inequality with (\ref{normeinfinie}). \fdem

\subsection{On the local expansion of the dominant eigenvalue $\lambda_t$  \label{local}}
 We explicit here the  local expansion  near $0$ of  the map $t \mapsto \lambda_t$:
  \begin{prop} \label{localexp}Under the hypotheses ${\bf H_0}{\bf \rm -}{\bf H_4}$,  there exists $C_\Gamma>0$ such that
 \begin{equation}\label{localexplambdat}
\lambda_t=1-C_\Gamma e^{-i\pi {\kappa\over 2}}t^\kappa L(1/t) \ (1+o(t)).
 \end{equation}
 \end{prop}
 \noindent {\bf Proof.} Recall first that $1$ is a simple eigenvalue of $\mathcal L_\delta$  with $\mathcal L_\delta h=h$ and  $\sigma(h)=1$; since $t\mapsto \mathcal L_{\delta+it} $ is continuous on $\mathbb R$, for $t$ closed to $0$ there exists a function $h_t\in \mathbb L_\omega$ such that $\mathcal L_{\delta+it}  h_t = \lambda_t h_t$, this function being unique if we impose the normalization condition $\sigma(h_t)=1$. The maps $t \mapsto \lambda_t$ and $t \mapsto h_t$ have the same type of continuity than $t \mapsto \mathcal L_{\delta+it} $ and one gets the identity
 $$
 \lambda_t=\sigma\left(\mathcal L_t h_t\right)=\sigma\left(\mathcal L_t h\right)
 +\sigma\bigl(    (\mathcal L_{\delta+it} -\mathcal L_\delta )(h_t-h))\bigr).$$
By the previous subsection, the second term  of this last  expression is $\preceq \Bigl(
t^\kappa L(1/t) \Bigr)^2$. It remains to precise the local behavior of the map $t\mapsto \sigma\left(\mathcal L_t h\right)$; one gets 
 $$\sigma\left(\mathcal L_t h\right)=1+ \sum_{l=0}^{r} \sigma_l$$ with 
$\displaystyle 
\sigma_0:= 
  \sum_{\beta \in \mathcal W_1}\int_{K_\beta}
h(\beta \cdot x)) e^{-\delta   {\mathfrak r}(\beta\cdot x)}(e^{-it    {\mathfrak r}(\beta\cdot x)}-1)\sigma(dx) 
 $
and, for $1\leq l\leq r$
$$
\sigma_l:= 
  \sum_{\beta \in \widehat{\mathcal P}_l}\int_{K_\beta}
h(\beta \cdot x)) e^{-\delta   {\mathfrak r}(\beta\cdot x)}(e^{-it    {\mathfrak r}(\beta\cdot x)}-1)\sigma(dx). 
$$
By (\ref{fourierpourbetK+1/K+L}), the terms $\sigma_l, l\neq r,$ are of the form $t^{\kappa }L(1/t) o(t).$
To control the   term  $\sigma_r$, one 
 sets $\Delta(n, x):= {\mathfrak r} (p_r^{n}\cdot x)-d({\bf x}_0, p_r^n\cdot {\bf x}_0)$ for any  $x \in  \overline{\partial X\setminus I_r} $ and $n \in \mathbb Z$. The following lemma readily implies that 
  the quantity $\Delta(n, x)$ tends to $-(x\vert  \xi_{r})_{\bf x_0}$ as $\vert n\vert \to +\infty$.

\begin{lem}
For any parabolic group  $\mathcal P:=<p>$ with   fixed point $\xi$, we have
$$
\mathcal B_x(p^{\pm n}\cdot  {\bf x}_0,  {\bf x}_0)=d_b ( {\bf x}_0,  p^n\cdot  {\bf x}_0)-2(\xi\vert  x)_ {{\bf x}_0}+\epsilon_x(n) 
$$
with $\displaystyle \lim_{n\to +\infty} \epsilon_x(n)=0$, the convergence being uniform on  compact sets of $\partial X\setminus\{\xi \}$.
\end{lem}

\noindent{\bf Proof.}   Let $({\bf x}_m)$ be a sequence of elements of $X$ converging to $x$.  We have
 \begin{eqnarray*}\mathcal B_x(p^{\pm n}\cdot  {\bf x}_0, {\bf x}_0)&=&\lim_m d (p^{\pm n}\cdot   {\bf x}_0, {\bf x}_m)-d (  {\bf x}_0, {\bf x}_m) \\
&=&  d (p^{\pm n}\cdot   {\bf x}_0,  {\bf x}_0)- \lim_m \Bigl( d ( {\bf x}_0, {\bf x}_m)+ d (p^{\pm n}\cdot  {\bf x}_0,  {\bf x}_0)-d (p^{\pm n}\cdot   {\bf x}_0, {\bf x}_m) \Bigr)\end{eqnarray*} 
with  $$ \lim_{n  } \Bigl(\lim_m  d ( {\bf x}_0, {\bf x}_m)+ d (p^{\pm n}\cdot  {\bf x}_0, {\bf x}_0)-d (p^{\pm n}\cdot  {\bf x}_0, {\bf x}_m) \Bigr)
= 2  \lim_{n  }   (  p^{\pm n}\cdot  {\bf x}_0 \vert  x)_{ {\bf x}_0}    =  2( \xi_h^-\vert  x)_{ {\bf x}_0}$$ and the conclusion follows as the Gromov product $( p^{\pm n} \cdot {\bf x}_0 \vert x )_{{\bf x}_0} $ tends uniformly to 
 $( \xi\vert  x)_{ {\bf x}_0}  $ on compacts of $\mathbb  S^1 \setminus  \{ \xi\}$.
\hfill{ $\Box$}

We write
 \begin{eqnarray*}
 \sigma_r&  = & \sum_{\vert n\vert \geq 2}
 e^{-\delta  d({\bf x}_0, p_r^{n}\cdot {\bf x}_0)}
\left(e^{-itd({\bf x}_0, p_r^{n}\cdot {\bf x}_0)}-1\right)\int_{\partial X\setminus I_r}
h(p_r^n \cdot x)) 
 e^{-\delta  \Delta (n, x)}
\sigma(dx)\\
 & &\quad + 
  \sum_{\vert n\vert \geq 2}
  e^{-\delta  d({\bf x}_0, p_r^n\cdot {\bf x}_0)}
 \int_{\partial X\setminus I_r}
h(p_r^n \cdot x)) 
 e^{-\delta  \Delta(n, x)}
 \Bigl( 
 e^{-it{\mathfrak r} (p_r^{n}\cdot  x)}-e^{-itd({\bf x}_0, p_r^n \cdot {\bf x}_0)}
 \Bigr)
\sigma(dx)\\
 &=& \sigma_{r1}+\sigma_{r2}.
  \end{eqnarray*}
  One gets $\displaystyle 
\int_{\partial X\setminus I_r}
h(p_r^n \cdot x)) 
 e^{-\delta  \Delta (n, x)}
\sigma(dx) >0$ for any  $\vert n\vert \geq 2$; by  (\ref{localexpansion}), there  exists $c >0$ such that 
$$\sigma_{r1}= -c e^{-i\pi {\kappa\over 2}}\Gamma(1-\kappa) t^\kappa L\Bigl({1\over t}\Bigr) (1+o(t)).$$
On the other hand
$$
\vert \sigma_{r2}\vert 
\leq \vert t\vert\sum_{\vert n\vert \geq 2}
  e^{-\delta  d({\bf x}_0, p_r^n\cdot {\bf x}_0)}\int_{\partial X\setminus I_r}
h(p_r^n \cdot x)) 
 e^{-\delta  \Delta(n, x)}\vert {\mathfrak r} (p_r^n\cdot  x)-d({\bf x}_0, p_r^n\cdot{\bf x}_0)\vert \sigma(dx)=O(t).
$$
Equality (\ref{localexplambdat}) follows immediately. 
 \fdem

 \subsection{Renewal theory and proof of Theorem \ref{countingkappa}} \label{renewal}

 For technical reasons (see for instance \cite{E}) which will appear in the control of the term $\widetilde W_j^{(3)}(R, \psi)$, we need to symmetrize the quantities $W_j( R, \psi )$ and   $W_j(s, R, \psi )$ setting  
  $$
 \widetilde W_j( R, \psi):= \sum_{k\geq 0} 
 \sum_{y\in \partial X\slash T^k\cdot y= x_j} e^{-\delta S_k\mathfrak r(y)}\Bigl(\psi(S_{k}{\mathfrak r}(y)-R)+
 \psi(-S_{k}{\mathfrak r}(y)-R)\Bigr)
 $$
 and
 $$
 \widetilde W_j(s, R, \psi):= \sum_{k\geq 0} s^k
 \sum_{y\in \partial X\slash T^k\cdot y= x_j} e^{-\delta S_k\mathfrak r(y)}\Bigl(\psi(S_{k}{\mathfrak r}(y)-R)+
 \psi(-S_{k}{\mathfrak r}(y)-R)\Bigr).
 $$
Notice that, when $\psi$ is a continuous function with compact support in $\mathbb R^+$,  the terms $\psi(-S_{k}{\mathfrak r}(y)-R)$ of these sums vanish for $R $ large enough, so that $\widetilde W_j(  R, \psi)= W_j(  R, \psi)$ and  $\widetilde W_j(s, R, \psi)= W_j(s, R, \psi)$ in this case.
By    (\ref{Wfourier}),   for any $0<s<1$, one gets,

\begin{eqnarray*}
\widetilde W_j(s, R, \psi)&:=& \sum_{k\geq 0} s^k
 \sum_{y\in \partial X\slash T^k\cdot y= x_j} e^{-\delta S_k\mathfrak r(y)}\Bigl(\psi(S_{k}{\mathfrak r}(y)-R)+
 \psi(-S_{k}{\mathfrak r}(y)-R)\Bigr)\\
&=&
{1\over 2\pi} \int_{\mathbb R}e^{itR}\hat{\psi}(t)\Bigl((I-s\mathcal L_{\delta+it})^{-1}-(I-s\mathcal L_{\delta-it})^{-1}\Bigr)1(x_j)dt.
\end{eqnarray*}
Fix $\epsilon_0>0$  and let $\rho (t)$ be a symmetric $C^\infty$-function on $\mathbb R$ which is equal to $1$ on a $[-\epsilon_0, \epsilon_0]$ and which  vanishes outside $[-2\epsilon_0, 2\epsilon_0]$; one thus decomposes $\widetilde W_j(s, R, \psi)
$ as 
$$\widetilde W_j(s, R, \psi)=\widetilde W^{(1)}_j(s, R, \psi)+ \widetilde W^{(2)}_j(s, R, \psi)+\widetilde W^{(3)}_j(s, R, \psi)$$
 with
\begin{eqnarray*}
\widetilde W^{(1)}_j(s,  R, \psi)&=& {1\over 2\pi} \int_{\mathbb R}e^{itR}\hat{\psi}(t)(1-\rho(t))\left(I-r\mathcal L_{\delta+it}  \right)^{-1}1(x_j)dt \\
& \ &\qquad +  {1\over 2\pi} \int_{\mathbb R}e^{itR}\hat{\psi}(t)(1-\rho(t))\left(I-r\mathcal L_{\delta-it} \right)^{-1}1(x_j)dt,\\
  \widetilde W^{(2)}_j(s, R, \psi)&=&
{1\over 2\pi} \int_{\mathbb R}e^{itR}\hat{\psi}(t)\rho(t)\left(\left(I-r\mathcal L_{\delta+it}  \right)^{-1}1(x_j)-
{\sigma(\partial X\setminus I_j)h(x_j)\over 1-r\lambda_t} \right) dt\\
&\ & \qquad + {1\over 2\pi} \int_{\mathbb R}e^{itR}\hat{\psi}(t)\rho(t)\left(\left(I-r\mathcal L_{\delta-it} \right)^{-1}1(x_j)-
{\sigma(\partial X\setminus I_j)h(x_j)\over 1-r\lambda_{-t}} \right) dt
\\
{\rm and} \  \widetilde W^{(3)}_j(s, R, \psi)&=&
{\sigma(\partial X\setminus I_j)h(x_j)\over 2\pi} \int_{\mathbb R}e^{itR}\hat{\psi}(t)\rho(t) \left({1\over 1-r\lambda_t}-{1\over 1-r\lambda_{-t}}\right)dt.
\end{eqnarray*}
Using Proposition \ref{localexp}, letting $s\to 1$, one gets $$  \widetilde W_j (R, \psi)=   \widetilde W_j^{(1)} ( R, \psi)+\     \widetilde W_j^{(2)}  ( R, \psi)+     \widetilde W_j^{(3)}  ( R, \psi)$$
 with
\begin{eqnarray*}\label{integrale1}
   \widetilde W_j^{(1)} (R, \psi)&=& \lim_{s\nearrow1}  \widetilde W_j^{(1)} (s,  R, \psi)\\
&=& {1\over 2\pi} \int_{\mathbb R}e^{itR}\hat{\psi}(t)(1-\rho(t))\Bigl(\left(I-\mathcal L_{\delta+it}  \right)^{-1}1(x_j)+\left(I- \mathcal L_{\delta-it}  \right)^{-1}1(x_j)\Bigr)dt,\\
   \widetilde W_j^{(2)} (R, \psi)&=& \lim_{s\nearrow1}  W^{(2)}_j(s,   R, \psi)\\
&=&
{1\over 2\pi} \int_{\mathbb R}e^{itR}\hat{\psi}(t)\rho(t)\left(\left(I-\mathcal L_{\delta+it}  \right)^{-1}1(x_j)-
{\sigma(\partial X\setminus I_j)h(x_j)\over 1-\lambda_t} \right) dt\\
&\ & \qquad + {1\over 2\pi} \int_{\mathbb R}e^{itR}\hat{\psi}(t)\rho(t)\left(\left(I-\mathcal L_{\delta-it}  \right)^{-1}1(x_j)-
{\sigma(\partial X\setminus I_j)h(x_j)\over 1- \lambda_{-t}} \right) dt
\\
{\rm and} \     \widetilde W_j^{(3)} ( R, \psi)&=& \lim_{s\nearrow1}  \widetilde W_j^{(3)} ( s, R, \psi)\\
&=&{\sigma(\partial X\setminus I_j)h(x_j)\over 2\pi} \int_{\mathbb R}e^{itR}\hat{\psi}(t)\rho(t) {\rm Re}\left({1\over 1-\lambda_t}\right)dt.
\end{eqnarray*}
The functions $t\mapsto \lambda_t$ has the same regularity as $t\mapsto L_{\delta+it} $;   by lemma \ref{regularitepourperturbationoperateurs} one can thus check that 
$
\psi_1: t \mapsto  \hat{\psi}(t)(1-\rho(t))\Bigl(\left(I-\mathcal L_{\delta+it}  \right)^{-1}1(x_j)+\left(I- \mathcal L_{\delta-it}  \right)^{-1}1(x_j)\Bigr)$
and 
\begin{eqnarray*}\psi_2: t &\mapsto&  \hat{\psi}(t)\rho(t)\Bigl(\left(I-\mathcal L_{\delta+it}  \right)^{-1}1(x_j)-
{\sigma(\partial X\setminus I_j)h(x_j)\over 1-\lambda_t} \\
&  & \qquad \qquad \qquad \qquad \qquad \qquad 
 +  \left(I-\mathcal L_{\delta-it} \right)^{-1}1(x_j)-
{\sigma(\partial X\setminus I_j)h(x_j)\over 1- \lambda_{-t}} \Bigr)
\end{eqnarray*}
satisfy the inequality $\vert \psi_k(s)-\psi_k(t)\vert \preceq \vert s-t\vert^\kappa L\left({1\over \vert s-t\vert}\right), k=1, 2$. This yields some information on the speed of convergence to $0$ of  their Fourier transform: indeed, for any $\theta < \kappa $, there exists $C_\theta>0$ such that 
\begin{equation}\label{majorationdeW1etW2}
\left\vert     \widetilde  W_j^{(1)} ( R, \psi)\right\vert\leq {C_\theta \over R^\theta}\qquad{\rm and} 
\qquad\left\vert     \widetilde W_j^{(2)} (  R, \psi)\right\vert\leq {C_\theta \over R^\theta}  
\end{equation}
(which readily implies $\displaystyle \lim_{R\to +\infty}R^{1-\kappa}L(R)    \widetilde W_j^{(1)} ( R, \psi)=\lim_{R\to +\infty}R^{1-\kappa}L(R)    \widetilde W_j^{(2)} ( R, \psi)=0$ as soon as $\theta\in ]1-\kappa, \kappa[$, which is possible since $\kappa>1/2$).
\\
On the other hand, by section 5 in \cite{E}, one gets for $1/2 <\kappa<1$
$$
\lim_{R\to +\infty}R^{1-\kappa} L(R)  \widetilde W_j^{(3)}( R, \psi)=C_j  \hat{\psi}(0)=C_j \int_{\mathbb R} \psi(x) dx 
$$
with $C_j = \sigma(\partial X\setminus I_j)h(x_j){\sin \pi \kappa\over \pi} $; notice that the value $h(x_j)$ is uniquely determined by the normalization  $\sigma(h)=1$ (see \cite{La} Theorem 4 for a detailed argument). This achieves the proof of Theorem \ref{countingkappa}.\fdem

\vspace{2cm}
Fran\c{c}oise Dal'bo \\ 
IRMAR, Universit\'e de Rennes-I\@,
 Campus de Beaulieu, 35042 Rennes Cedex.\\ 
 mail : francoise.dalbo@univ-rennes1.fr
\vspace{5mm} 
 
  Marc Peign\'e \& Jean-Claude Picaud \\
LMPT, UMR 6083, Facult\'e des Sciences et Techniques, Parc de Grandmont, 37200 Tours. \\
mail : peigne@univ-tours.fr,\\
\hspace{5mm} jean-claude.picaud@univ-tours.fr  
\vspace{5mm}

 Andrea Sambusetti\\
        Istituto di Matematica G. Castelnuovo, Sapienza Universit\`a di Roma,      P.le Aldo Moro 5 - 00185 Roma. \\ 
        mail : sambuset@mat.uniroma1.it

\end{document}